\documentclass{birkjour_t2}
 \usepackage[utf8]{inputenc} 
 \usepackage{placeins} 
\usepackage{csquotes}
\usepackage{enumitem} 
\usepackage{xcolor} 
\usepackage{amsmath}
\usepackage{amssymb}
\usepackage{amsthm} 

\def\R{{\mathbb{R}}}
\def\N{{\mathbb{N}}}

\theoremstyle{plain}
\newtheorem{theo}{Theorem}

\newtheorem{prop}[theo]{Proposition}
\newtheorem{lemma}[theo]{Lemma}
\newtheorem{corol}[theo]{Corollary}

\theoremstyle{definition}

\theoremstyle{remark}
\newtheorem{rmk}[theo]{Remark}

\usepackage[]{graphicx}

\usepackage{float}
\usepackage[labelfont=bf]{caption}
\usepackage{subcaption}
\def\interior{{\rm int}\,}

\def\ninfty{{n \rightarrow +\infty}}

\newcommand{\normX}[1]{\|{#1}\|_X}

\newcommand{\normLp}[2]{\|{#1}\|_{L^{#2}(0,1)}}

\newcommand{\infimum}[1]{\inf\limits_{#1}}
\newcommand{\supremum}[1]{\sup\limits_{#1}}
\newcommand{\maximum}[1]{\max\limits_{#1}}
\newcommand{\limit}[1]{\lim\limits_{#1}}
\newcommand{\limesinferior}[1]{\liminf\limits_{#1}}

\newcommand{\toinfty}[1]{{#1} \rightarrow +\infty}

\def\BVP{{(\ref{BVP_merge_eq})--(\ref{BVP_merge_cond_II})}}

\def\mezeraA{\vspace{4 ex}}

\usepackage{hyperref}

\begin{document}

\title{Variational and numerical aspects of a system of ODEs with concave-convex nonlinerities}

\author{Oscar Agudelo}
\email{oiagudel@ntis.zcu.cz}
 \address{Department of Mathematics and NTIS, Faculty of Applied 
Sciences, University of~West Bohemia, Univerzitn\'{\i} 8, 301~00~Plze\v{n}, Czech Republic, ORCID: 0000-0002-2588-9999}

\author{Gabriela Holubov\'{a}}
\email{gabriela@kma.zcu.cz}
\address{Department of Mathematics and NTIS, Faculty of Applied Sciences, University of~West Bohemia, Univerzitn\'{\i} 8, 301~00~Plze\v{n}, Czech Republic, ORCID: 0000-0003-1127-3381}

\author{Martin Kudl\'{a}\v{c}}
\email{kudlacm@kma.zcu.cz}
\address{Department of Mathematics and NTIS, Faculty of Applied Sciences, University of~West Bohemia, Univerzitn\'{\i} 8, 301~00~Plze\v{n}, Czech Republic, ORCID: 009-0007-1749-3314}

\maketitle

{\bf Abstract.} In this work we discuss a Hamiltonian system of ordinary differential equations under Dirichlet boundary conditions. The system of equations in consideration features a mixed (concave-convex) power nonlinearity depending on a positive parameter $\lambda$. We show multiplicity of nonnegative solutions of the system for a certain range of the parameter  $\lambda$ and we also discuss regularity and symmetry of nonnegative solutions of the system. Besides, we present a numerical strategy aiming at the exploration of the optimal range of $\lambda$ for which multiplicity of solutions holds. The numerical experiments are based on the {Poincar\'{e}-Miranda} theorem and the {shooting method}, which have been lesser explored for systems of ODEs. Our work is motivated by the works of Ambrosetti et al in \cite{ambrosettibrezis} and Moreira dos Santos in \cite{dossantos}. 

\vspace{10 pt}


\vspace{10 pt}

\section{Introduction}\label{sec:intro}

In this work we study the system of ordinary differential equations (ODEs for short)
\begin{equation}
		\left\{ 
			\begin{aligned}	
						-u'' &= \lambda\, |v|^{r-1}v + |v|^{p-1}v \quad &\hbox{in} \quad &(0,1),\\ 
						-v'' &= |u|^{q-1} u \quad &\hbox{in} \quad &(0,1)
			\end{aligned} 
		\right.
		\label{BVP}
	\end{equation}
with \emph{Dirichlet boundary conditions} 
\begin{equation}
\label{eqn:BVP_Dirch}
u(0) = u(1)=0 \quad \hbox{and} \quad v(0)= v(1) = 0.
\end{equation}

We consider $\lambda$ to be a positive parameter and the exponents $p,q$ and $r$ being such that
\begin{equation}\label{hyp:pqr}
0<r<\frac{1}{q} \quad \mbox{and} \quad p > \max\left\{1,\frac{1}{q}\right\}. 
\end{equation}

Systems of equations as in \eqref{BVP} appear naturally in the study of population dynamics, fluid dynamics and stellar structure in astrophysics. Concerned with the latter, related systems of equations have been explored since the ninetieth century to study for instance the density of a gas sphere (see e.g. \cite{chandrasekhar,kippenhaln}).

Notice that system \eqref{BVP} features concave-convex polynomial nonlinearities. When $0<r<1<p$, the nonlinearity in the first equation of \eqref{BVP}
is concave near the origin and convex at infinity. When $0<q<1<r$, the nonlinearity in the first equation of \eqref{BVP} stays convex while concavity appears in the second equation of \eqref{BVP}.  

In \cite{ambrosettibrezis}, the authors study existence, nonexistence and multiplicity of nonnegative solutions of
the single equation
\begin{equation}
-v'' = \lambda\, |v|^{r-1} v + |v|^{p-1} \, v \quad \hbox{in} \quad (0,1)		\label{eqn:BVP_ABC}
	\end{equation}
with {Dirichlet boundary conditions} 
\begin{equation}
\label{eqn:BVP_Dirch_ABC}
v(0)= v(1) = 0
\end{equation}
and with respect to the nonnegative parameter $\lambda$, 
assuming that $0<r<1<p<+\infty$. Actually, the study in \cite{ambrosettibrezis} treats mainly the higher dimensional case, where technical issues related with the non-compactness of certain Sobolev embeddings arise. In \cite{dossantos}, the author studies system \eqref{BVP}-\eqref{eqn:BVP_Dirch} in higher dimensions. The spirit of the results in \cite{dossantos} is similar to the one in \cite{ambrosettibrezis}. In \cite{dossantos} the range of values $\lambda$ for which existence and multiplicity of solutions of \eqref{BVP}-\eqref{eqn:BVP_Dirch} are studied is not optimal and we remark that the techniques in \cite{dossantos} differ from the ones in \cite{ambrosettibrezis}. In \cite{agudelo2023hamiltonian}, the authors discuss the existence of minimal solutions for a related Hamiltonian system of equations. Nonexistence of solutions is also discussed. We refer the reader to the surveys \cite{de2008semilinear} and \cite{ruf2008superlinear} and references therein for a description of some of the known results related to systems of differential equations and the techniques used to treat them.

\vskip 3pt	
Motivated mainly by \cite{ambrosettibrezis} and \cite{dossantos}, we explore in dimension one the \emph{existence} and \emph{multiplicity} of \emph{nonnegative classical solutions} of \eqref{BVP}-\eqref{eqn:BVP_Dirch} with respect to the parameter $\lambda$. By a classical nonnegative  solution to \eqref{BVP}-\eqref{eqn:BVP_Dirch} we mean a pair $(u,v)$ of functions such that $u,v\in C^2(0,1)\cap C[0,1]$, $u,v\geq 0$ in $(0,1)$ and \eqref{BVP} and \eqref{eqn:BVP_Dirch} being satisfied point-wise.

\vskip 3pt
Our main result is the following. 

\begin{theo}\label{theo:main_theo_1}		There exists a positive constant $\lambda_0$ such that for any $\lambda \in (0, \lambda_0)$, the system \eqref{BVP}-\eqref{eqn:BVP_Dirch} possesses at least
two distinct nontrivial nonnegative classical solutions.
\end{theo}

We now explain briefly the strategy of the proof of Theorem \ref{theo:main_theo_1}. Based on the method of reduction by inversion (see e.g. \cite{clement1997homoclinic,hulshof1996asymptotic}), and proceeding as in \cite{dossantos}, we reformulate \eqref{BVP}-\eqref{eqn:BVP_Dirch} as a single fourth order ODE with Navier boundary conditions (see BVP \eqref{BVP_merge_eq}-\eqref{BVP_merge_cond_II}). This fourth order BVP has variational structure (for the corresponding energy functional see \eqref{energy_BVP}). Following partly the approach in \cite{ambrosettibrezis}, we apply the \emph{Ambrosetti-Rabinowitz Mountain Pass Theorem} and a \emph{local minimization argument} to the associated energy functional to prove the existence of two distinct weak solutions for sufficiently small $\lambda$ (see Propositions \ref{TH: MP solution} and \ref{TH: solution near the origin}, respectively).

We remark that in contrast with the fibering method used in \cite{dossantos}, this variational approach easily adapts to study existence of nonnegative solutions of BVP of the type \eqref{BVP}-\eqref{eqn:BVP_Dirch} with more general nonlinearities. For instance, one may consider nonlinearities with similar mixed polynomial growth, but such that this feature does not come from autonomous terms. For the sake of clarity and emphasis of the main ideas in our exposition, we work directly with the BVP \eqref{BVP}-\eqref{eqn:BVP_Dirch}.

\vskip 3pt
In addition to discussing existence and multiplicity, we prove the regularity of solutions (see Proposition \ref{theo:regl_4th_order_BVP}) and provide a lower estimate of the optimal value of $\lambda_0$. We finish our theoretical discussion by using the method of moving planes to prove symmetry of the solutions (see Proposition \ref{theo:symmetry_slns}). 

\vskip 3pt
 Although we use quite standard functional analytical tools, the nonlinear terms in \eqref{BVP} require a fine and careful treatment. Due to the techniques and restrictions of Sobolev embeddings used in \cite{dossantos}, the results concerned with existence of solutions for the higher-dimensional version of system \eqref{BVP} do not directly translate into the one dimensional setting. Nonetheless, as in \cite{dossantos}, the variational techniques used in this work do not provide the optimal quantitative information concerning the range of values 
$\lambda$ for which existence and multiplicity of nonnegative solutions hold true.
Its analytic description (depending on parameters $p$, $q$, $r$) remains to be open.

\vskip 3pt
We therefore further explore this matter further using a numerical implementation motivated by a combination of the standard \emph{Shooting method} and the heuristics of the \emph{Poincar\'{e}-Miranda Theorem}. To the best of our knowledge, this strategy has not been explored yet for systems of the type \eqref{BVP}-\eqref{eqn:BVP_Dirch} and hence we believe the numerical illustration presented in Section \ref{sec:numerical_exp} is a novel, insightful and versatile approach that easily adapts to more general systems of ODEs. 

\begin{figure}[htb]
		\centering
		\begin{picture}(300,210)
			\put(0,0){\includegraphics[width = 10 cm]{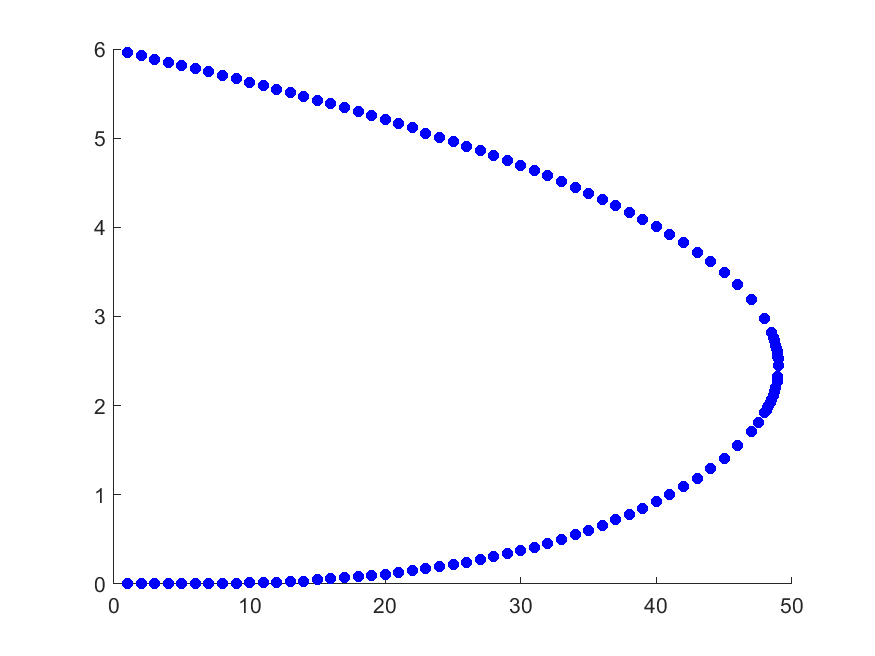}}
			\put(150,0){\small{$\lambda$}}
			\put(0,90){\rotatebox{90}{\small{$\maximum{x\in[0,1]}v(x)$}}}
		\end{picture}
		\caption{Bifurcation diagram for the parameter $\lambda$. Here, $p=3$, $q=3/2$, $r=1/3$}
		\label{fig:bifurcation diagram}
	\end{figure}

In particular, in line with the statement of Theorem \ref{theo:main_theo_1}, 
we obtain
a bifurcation diagram showcasing the dependence of the $L^{\infty}(0,1)$-norm of the 
$v$-component of a solution of \eqref{BVP}-\eqref{eqn:BVP_Dirch} with respect to the parameter $\lambda$ (see Fig. \ref{fig:bifurcation diagram}). 	

Borrowing some terminology from the Bifurcation Theory, our results suggest the following: there exists $\lambda_{{\rm bif}}>0$ such that for $\lambda \in [0,\lambda_{{\rm bif}})$ there are at least two solutions -- a stable one $v_{\lambda}$ and the unstable one $v^{\lambda}$. At the point of bifurcation $\lambda=\lambda_{{\rm bif}}$, the solutions $v_{\lambda}$ and $v^{\lambda}$ coincide. We remark that, as shown in Theorem 1.3 in \cite{agudelo2023hamiltonian}, when $0<r<\frac{1}{q}<1<p$, there exists $\Lambda_0\in (0,\infty)$ such that for any $\lambda > \Lambda_0$ the system \eqref{BVP} has no nontrivial nonnegative solutions. It is expected that the smallest of such $\Lambda_0$ coincides with $\lambda_{{\rm bif}}$.\\

The paper is organized as follows. Section \ref{section: preliminaries} introduces preliminary results and notions required in subsequent sections. In Section \ref{sec:setting} we set the functional analytic framework to study the problem \eqref{BVP}-\eqref{eqn:BVP_Dirch}. This section also discusses the regularity of solutions. To carry out the proof of the main result, Sections \ref{sec:first_sln} and \ref{sec:second_sln} provide existence of two weak solutions of the fourth order BVP (see \eqref{BVP_merge_eq}-\eqref{BVP_merge_cond_II}) via the Mountain Pass Theorem and a local minimization argument for the associated energy functional (see \eqref{energy_BVP}), respectively. Section \ref{sect:PS} deals with compactness properties of approximating sequences of solutions, namely, the Palais-Smale condition.
The proof of Theorem \ref{theo:main_theo_1} is carried out in Section \ref{sec:proof_of_main_theorem} and Section \ref{section:Symmetry} is devoted to the proof of Proposition \ref{theo:symmetry_slns} describing symmetry of the solutions. The last section presents various numerical illustrations with detailed description of our numerical strategy.

\section{Preliminaries} \label{section: preliminaries}

Let $s\in [1,\infty]$ and let $I$ be a bounded open interval. In what follows, $L^s(I)$ denotes the Lebesgue space of measurable functions $w:I\to \R$ endowed with the norm 
$$
\|w\|_{L^s(I)}:=\left\{\begin{aligned}
	\left( \int\limits_I |w|^s \mbox{d}s \right)^{\frac{1}{s}}&, \qquad s\in [1,+\infty),\\
	{\rm ess}\sup \limits_{x\in I}|w(x)|&,\qquad s=+\infty.
\end{aligned}
\right.
$$ 

When $I=(0,1)$ and for $k\in \mathbb{N}$, we also consider the Sobolev space $W^{k,s}(0,1)$, which consists on all functions $u \in L^s(0,1)$ such that for any $i \in \{1,\ldots,k\}$ the $i-$th weak derivative of $u$, $u^{(i)}$, belongs to $ L^s(0,1)$. The space $W^{k,s}(0,1)$ is endowed with the norm \begin{equation*}
	\|u\|_{W^{k,s}(0,1)}:= \|u\|_{L^s(0,1)} + \sum\limits_{i = 1}^{k} \|u^{(i)}\|_{L^s(0,1)}.
\end{equation*}
Recall that the Sobolev space $W^{k,s}(0,1)$ is a Banach space. Furthermore, $W^{k,s}(0,1)$ is separable for $s \in [1, \infty)$ and reflexive for $s \in (1, \infty)$ (see Theorems 3.3 and 3.6 in \cite{adams}, p. 60--61).

\vspace{10 pt}

Let $s>1$ and let $W_0^{1,s}(0,1)$ denote the closure of $C_c^1(0,1)$ in $W^{1,s}(0,1)$ with respect to $\|.\|_{W^{1,s}(0,1)}$. A convenient description of the space $W^{1,s}_0(0,1)$ is the following (see \cite{brezis}, Th. 8.12, p. 217): for any $u \in W^{1,p}(0,1)$,  
\vspace{2 pt} 
\begin{center}$u\in W^{1,p}_0(0,1)$ if and only if $u = 0$ on $\partial I$.\end{center}

\begin{lemma}[Morrey's inequality revisited] \label{lemma: iterated Morrey's}
	If $u\in W^{2,s}(0,1) \cap W^{1,s}_0(0,1)$, then $$\|u\|_{L^\infty(0,1)} \leq \frac{1}{2} \|u''\|_{L^s(0,1)}.$$
\end{lemma}
\begin{proof}
	First, notice that if $w \in L^\infty(0,1)$, then for any $x,y \in (0,1)$ with $x<y$, we estimate
	\begin{equation}
		\|w\|_{L^s(x,y)} = \left( \int\limits_x^y |w|^s \mbox{d}s \right)^{\frac{1}{s}} \leq (y-x)^{\frac{1}{s}}\|w\|_{L^\infty(x,y)}
		\label{equation: Ls norm <= Linfty norm}
	\end{equation}

	Also, if $w \in W^{1,s}(0,1)$, then for any $x,y\in [0,1]$ with $x<y$, $w\in W^{1,s}(x,y)$. Even more, from the Morrey's inequality (\cite{evans}, Th. 4, p. 280),
	\begin{equation}
	|w(y)-w(x)| \leq  \, (y-x)^{1 -\frac{1}{s}} \|w'\|_{L^s(x,y)}.\label{Linfty norm of u <= Lp norm of u'}
	\end{equation}

	Finally, if $w \in W_0^{1,\sigma}(0,1)$ with $\sigma \geq 2$, then, for any $x \in (0,1)$, \eqref{Linfty norm of u <= Lp norm of u'} implies
	$$
	|w(x)|^{\sigma} \leq \, x^{\sigma - 1} \|w'\|_{L^{\sigma}(0,x)}^{\sigma} \quad \mbox{and} \quad
	|w(x)|^{\sigma} \leq \, (1-x)^{\sigma - 1} \|w'\|_{L^{\sigma}(x,1)}^{\sigma}.
	$$
	Multiplying the left inequality by $(1-x)^{\sigma -1}$ and the right inequality by $x^{\sigma -1}$ and summing up, yields
	$$
	\left((1-x)^{\sigma - 1} + x^{\sigma -1} \right) |w(x)|^{\sigma} \leq \, x^{\sigma - 1}(1-x)^{\sigma - 1}  \|w'\|_{L^{\sigma}(0,1)}^{\sigma}.
	$$
	Applying Jensen and A-G inequalities, we obtain
	\begin{equation}
		|w(x)| \leq  \, \frac{1}{2} \|w'\|_{L^{\sigma}(0,1)}.\label{Linfty norm of u <= Lp norm of u' for W0}
	\end{equation}
	
    \mezeraA
    
	Now, let $u\in W^{2,s}(0,1)\cap W^{1,s}_0(0,1)$, $s > 1$. Then $u \in W^{1,\sigma}_0(0,1)$ with any $\sigma \geq 2$ and thus \eqref{Linfty norm of u <= Lp norm of u'}, \eqref{Linfty norm of u <= Lp norm of u' for W0}, and Rolle's Theorem (see \cite{zorich}, p. 215) yield
	$$\|u\|_{L^\infty(0,1)} \leq \frac{1}{2} \|u'\|_{L^{\sigma}(0,1)} \leq \frac{1}{2} \|u'\|_{L^{\infty}(0,1)} \leq \frac{1}{2} \|u''\|_{L^s(0,1)},$$ in other words the statement of this lemma.
\end{proof}

Next, let us fix $q>0$  and introduce
\begin{equation}
X:=W^{2,\frac{q+1}{q}}(0,1)\cap W^{1,\frac{q+1}{q}}_0(0,1).
 \label{definition of X gamma}
\end{equation}

Observe that $u \in X$ implies that
$u \in C^{1,\frac{1}{q+1}}[0,1]$, $u''\in L^{\frac{q+1}{q}}(0,1)$ and $u(0) = u(1) = 0$.
 
It is readily verified using Lemma \ref{lemma: iterated Morrey's} that the mapping
\begin{equation*}
    \begin{aligned}
        \normX{v} := \left( \int\limits_0^1 |v''(x)|^\frac{q+1}{q} \, \mbox{d}x \right)^\frac{1}{q+1}
    \end{aligned}
\end{equation*}
is a norm in $X$. From now on (unless stated otherwise), we endow $X$ with this norm.

\begin{rmk}     \label{lemma: L infinity norm <= X norm of u} Notice that Lemma \ref{lemma: iterated Morrey's} implies the following. If $u \in X$, then for any $s \in (1, \infty]$, $u \in L^s(0,1)$ and
    \begin{equation}
         \normLp{u}{s} \leq\normLp{u}{\infty} \leq \frac{1}{2} \normX{u}.
        \label{EQ: infty norm of u <= X norm of u}
    \end{equation}
    
    Moreover, 
    the linear embedding
    \begin{equation}
    i: X \rightarrow L^s(0,1), \qquad 
  i(u) := u
        \label{injection from X to Ls}
    \end{equation}
    is continuous in $X$ and $\|i\| \leq \frac{1}{2}$.
\end{rmk}

\begin{lemma} \label{lemma: norm equivalence, Banach, Reflexive}
    The normed space $X$ has the following properties:
    \begin{enumerate}[label=(\alph*)]
        \item the norms $\normX{.}$ and $\|.\|_{W^{2,\frac{q+1}{q}}(0,1)}$ are equivalent in $X$,\label{enum a}
        \item space $X$ is Banach,\label{enum b}
        \item the space $X$ is reflexive, \label{enum c}
        \item the space $X$ is compactly embedded into $C^1[0,1]$. \label{enum d}
    \end{enumerate}
\end{lemma}
Claims \ref{enum a}-\ref{enum c} are thoroughly proved in \cite{kudlacDiploma}, the last property of $X$ is obtained by iterating of compact embedding of $W^{1,\frac{q+1}{q}}(0,1)$ into $C[0,1]$ (see \cite{brezis}, Th. 8.2, p. 204).

\section{Functional analytic setting} 
\label{sec:setting}
In this section we discuss several formulations of 
\eqref{BVP}-\eqref{eqn:BVP_Dirch} and concepts of its solutions.
Recall that $p$, $q$ and $r$ satisfy \eqref{hyp:pqr}. \\ 

First of all, in order to capture nonnegative solutions, instead of working directly with \eqref{BVP}-\eqref{eqn:BVP_Dirch}, we consider the system	
    \begin{equation}
		\left\{ 
		\begin{aligned}	
			-u'' &= \lambda\, (v_+)^r + (v_+)^p \quad &\hbox{in} \quad &(0,1),\\
			-v'' &= |u|^{q-1} \, u \quad &\hbox{in} \quad &(0,1),\\
			u(0) &= u(1) = v(0) = v(1) = 0,
		\end{aligned} 
		\right.
		\label{BVP positive parts}
	\end{equation}
	where $t_+:=\max\{t,0\}$ for $t\in \R$.
Indeed, let $u, v\in C^2(0,1)\cap C[0,1]$ be such that the pair $(u,v)$ is a classical solutions of \eqref{BVP positive parts}, that is, all the identities in \eqref{BVP positive parts} hold in point-wise sense. First notice that if the pair $(u,v)$ is not the trivial vector function $(0,0)$, then both $u$ and $v$ are nontrivial. Next, from the first equation in \eqref{BVP positive parts}, $u$ is concave in $[0,1]$. Moreover, since $u(0) = u(1) = 0$, $u$ is nonnegative in $[0,1]$. Arguing in a similar manner with the second equation, using that $u\geq 0$ in $(0,1)$, the same holds for $v$. In conclusion, any classical nontrivial solution $(u,v)$ of (\ref{BVP positive parts}) is such that $u$ and $v$ are positive in $(0,1)$ and hence it is also a solution of \eqref{BVP}-\eqref{eqn:BVP_Dirch}. Clearly, any classical solution $(u,v)$ of (\ref{BVP}) with $u\geq 0$ and $v\geq 0$ in $[0,1]$ solves also (\ref{BVP positive parts}).\\
 
Also, from the second equation in (\ref{BVP positive parts}) we have 
	\begin{equation}
		u = -|v''|^{\frac{1}{q}-1} \, v'' \in C^2(0,1).
	\label{u vyjadreno pomoci v}
	\end{equation}
	
Plugging (\ref{u vyjadreno pomoci v}) into the first equation in (\ref{BVP positive parts}), system \eqref{BVP}-\eqref{eqn:BVP_Dirch} reduces to
the \emph{fourth order} equation
\begin{equation}
		\frac{\mbox{d}{ }^2}{\mbox{d}{x^2}} \left(|v''|^{\frac{1}{q}-1} v'' \right) = \lambda \, v_+^{r} + v_+^{p} \quad \hbox{in}\quad (0,1)
	\label{BVP_merge_eq}
	\end{equation}
with the \emph{Navier boundary conditions} 
	\begin{equation}
		v(0) = v(1) = 0 \quad \hbox{and}  \quad v''(0) = v''(1)=0.\label{BVP_merge_cond_II}
	\end{equation}

A \emph{classical solution} of \eqref{BVP_merge_eq}-\eqref{BVP_merge_cond_II} is a function $v\in C^2[0,1]$ such that $|v''|^{\frac{1}{q}-1}v''\in C^2(0,1)$ and \eqref{BVP_merge_eq}-\eqref{BVP_merge_cond_II} are satisfied point-wise. Observe that nonnegative classical solutions of the system \eqref{BVP}-\eqref{eqn:BVP_Dirch} are in correspondence with classical solutions of the BVP \eqref{BVP_merge_eq}-\eqref{BVP_merge_cond_II}.\\

Besides the classical setting for solutions, we will also use the concept of weak solutions. We consider the space $X$ defined in \eqref{definition of X gamma}.
By a \emph{weak solution} to \eqref{BVP_merge_eq}-\eqref{BVP_merge_cond_II} we understand a function $v\in X$ such that for any $\varphi\in X$,\begin{equation}
\int\limits_0^1 |v''|^{\frac{1}{q}-1} \, v'' \, \varphi'' \mbox{d}x = \int\limits_0^1 \left( \lambda \, v_+^{r} \, +\, v_+^{p} \right)\varphi\,\mbox{d}x.
\label{weak formulation}
\end{equation}
In the following statements we show that weak and classical solutions of the BVP \eqref{BVP_merge_eq}-\eqref{BVP_merge_cond_II} coincide.
 
\begin{prop}\label{theo:regl_4th_order_BVP}
		If $v \in X$ is a weak solution of \eqref{BVP_merge_eq}-\eqref{BVP_merge_cond_II}, then $v$ is also a classical solution of \eqref{BVP_merge_eq}-\eqref{BVP_merge_cond_II}.	
\end{prop}

\begin{proof}
Let $v\in X$ be a weak solution of \BVP. Write $w:=|v''|^{\frac{1}{q}-1}v''$ and $h:=\lambda v_+^{r} + v_+^{p} $ a.e. in $(0,1)$. Observe that $w\in L^{q+1}(0,1)$, $h\in C[0,1]$ and from \eqref{weak formulation} for any $\varphi \in X$
$$
\int_0^1 w\,\varphi'' \mbox{d}x =\int_{0}^1 h\,\varphi \,\mbox{d}x.
$$

Let $\omega\in C^2[0,1]$ solve the  BVP
$$
\omega''=h \quad \hbox{in}\quad (0,1), \qquad \omega(0)=\omega(1)=0. 
$$ 

We prove that $w=\omega$ a.e. in $(0,1)$. Once this is proven, we would conclude that $|v''|^{\frac{1}{q}-1}v'' \in C^2[0,1]$.

\vskip 3pt
To prove the claim, write $\vartheta:= |w-\omega|^{q-1}(w-\omega)$ a.e. in $(0,1)$. Since $w\in L^{q+1}(0,1)$, we have $\vartheta \in L^{\frac{q+1}{q}}(0,1)$. Let $\varphi \in W^{2,\frac{q+1}{q}}(0,1)$ solve  the BVP
$$
\varphi''=\vartheta \quad \hbox{in}\quad (0,1), \qquad \varphi(0)=\varphi(1)=0.
$$

Observe that $\varphi \in C^1[0,1]$ and from the boundary conditions for $\varphi$ we find that $\varphi \in X$. Using $\varphi$ as a test function in \eqref{weak formulation}, we calculate
$$
\begin{aligned}
\int_0^1|w-\omega|^{q+1}\mbox{d}x=&\int_{0}^1 (w-\omega)\vartheta\,\mbox{d}x \qquad\big(\hbox{definition of }\vartheta\big)\\
=&\int_0^1(w-\omega)\varphi''\mbox{d}x\qquad\big(\hbox{choice of }\varphi\big)\\
=&\int_0^1 w\,\varphi''\mbox{d}x - \int_0^1 \omega''\varphi \,\mbox{d}x  \qquad\big(\hbox{integration by parts}\big)\\
=&\int_0^1 h\,\varphi \,\mbox{d}x - \int_0^1 h \,\varphi \,\mbox{d}x 
= 0.
\end{aligned}
$$

Thus, $w=\omega$ a.e. in $(0,1)$ and this proves the claim. Since $v''=|w|^{q-1}w$ in $[0,1]$, then $v''\in C[0,1]$. Consequently, $v\in C^2[0,1]$. This completes the proof of the proposition. 
\end{proof}

For an alternative proof of Proposition \ref{theo:regl_4th_order_BVP} (under slightly stronger assumptions on the exponent $q$, motivated by the arguments presented in \cite{fucik}), we refer the reader to \cite{kudlacDiploma}. The following corollary is a direct consequence of Proposition \ref{theo:regl_4th_order_BVP}.

\vskip 3pt

\begin{corol}\label{coro:system_to_4thorder}
Let $v\in X\setminus\{0\}$ be a nonnegative (weak) solution of \eqref{BVP_merge_eq}-\eqref{BVP_merge_cond_II} and set $u=-|v''|^{\frac{1}{q}-1}v''$ a.e. in $(0,1)$. Then $u,v\in C^2[0,1]$, $u,v>0$ in $(0,1)$ and the pair $(u,v)$ is a classical solution of \eqref{BVP}-\eqref{eqn:BVP_Dirch}. 
\end{corol}

Summarizing the above results, classical nonnegative nontrivial solutions to \eqref{BVP}-\eqref{eqn:BVP_Dirch} are in correspondence with the 
nontrivial weak solutions of \eqref{BVP_merge_eq}-\eqref{BVP_merge_cond_II}. As we have already announced, we will use a variational approach to find them.\\

Let us consider the energy functional $J: X \rightarrow \R$ defined by 
\begin{equation}
J(v) := \frac{q}{q+1} \, \int\limits_0^1 |v''|^{\frac{q+1}{q}} \, \mbox{d}x - \frac{1}{r+1} \, \int\limits_0^1 \lambda \, v_+^{r+1}  \, \mbox{d}x - \frac{1}{p+1} \, \int\limits_0^1 v_+^{p+1} \, \mbox{d}x.
\label{energy_BVP}
\end{equation}
	
	From Remark \ref{lemma: L infinity norm <= X norm of u}, $J$ is well-defined. It is also standard to verify that $J\in C^1(X)$ with \begin{equation}
		DJ(v)\varphi = \int\limits_0^1 |v''|^{\frac{1}{q}-1}v'' \varphi'' \mbox{d}x - \int\limits_0^1 \left( \lambda v_+^{r}   + v_+^{p}\right)\varphi \, \mbox{d}x 
        \quad \hbox{for }v,\varphi\in X.
		\label{EQ: Gateaux derivative of J}
	\end{equation}

	Thus, $DJ(v)= 0$ if and only if $v$ is a weak solution of \BVP. We refer the reader to \cite{kudlacDiploma} for details. 
	
	For convenience, we write $J(v)$ as
	\begin{equation}
		J(v) = \frac{q}{q+1}\normX{v}^{\frac{q+1}{q}} - \frac{\lambda}{r+1} \normLp{v_+}{r+1}^{r+1} - \frac{1}{p+1} \normLp{v_+}{p+1}^{p+1}.
		\label{J(v) norms}
	\end{equation}
In the following sections, we examine the critical points of $J$.

\section{Mountain pass solution}\label{sec:first_sln}
  
In this part, we use the \emph{Mountain Pass Theorem} to show that for $\lambda>0$ small, there exists at least one solution of \BVP. 
In what follows, $X^*$ denotes the topological dual space of $X$ with the topology induced by the norm in $X$. Recall that $p,q,r$ satisfy \eqref{hyp:pqr}. \\

First, we remark that the functional $J$ satisfies the Palais-Smale condition  (PS-condition for short) in $X$, i.e., for any $c\in \R$ and for any sequence ${(u_n)} \subset X$ such that
	\begin{equation}
		J(u_n) \rightarrow c, \quad \hbox{and}\quad \|DJ(u_n)\|_{X^*} \rightarrow 0,
		\label{Assumption Palais-Smale}
	\end{equation}
there exists a subsequence $(u_{n_k})\subset (u_n)$ that converges strongly in $X$ (see e.g. \cite{drabekmilota}). The detailed proof of this fact is postponed until Section \ref{sect:PS}.\\

Next, we verify that the functional $J$ has the Mountain Pass geometry. That is the scope of the following lemmas.

\begin{lemma} \label{lemma: J(v) positive on a sphere}
	There exist positive constants $T$ and $\lambda_0 = \lambda_0(T,p,q,r)$ such that for all $\lambda \in (0,\lambda_0)$ there exists $C = C(T,p,q,r,\lambda)$ so that for any 
	$v \in X$ satisfying $\|v\|_X = T$, $J(v) \geq C > 0$.
\end{lemma}
\begin{proof}
	Observe that for any $v \in X$, $v_+\leq |v|$ in $(0,1)$. Using \eqref{J(v) norms} and \eqref{EQ: infty norm of u <= X norm of u} from Remark~\ref{lemma: L infinity norm <= X norm of u} we estimate
 \begin{equation}
	\begin{aligned}
  		J(v) & \geq \frac{q}{q+1}\normX{v}^{\frac{q+1}{q}} - \frac{\lambda}{2^{r+1}(r+1)} \normX{v}^{r+1} - \frac{1}{2^{p+1}(p+1)} \normX{v}^{p+1}\\
		& \geq \normX{v}^{r+1} \left( \frac{q}{q+1} \, \normX{v}^{\frac{1}{q}-r} - \frac{\lambda}{2^{r+1}(r+1)} - \frac{1}{2^{p+1}(p+1)} \, \normX{v}^{p-r} \right).	    
	\end{aligned}
  \label{J(v) estimated using advanced Poincare}
	\end{equation}

	Consider the function
	$$h(t) := \frac{q}{q+1} \, t^{\frac{1}{q}-r} - \frac{1}{2^{p+1}(p+1)} \, t^{p-r}, \quad t \geq 0.$$

	A direct calculation yields that 
	$$h'(t) = t^{\frac{1}{q}-r-1} \left( \frac{1-qr}{q+1} - \frac{p-r}{2^{p+1}(p+1)} \, t^{p-\frac{1}{q}} \right) \quad \hbox{for }t > 0.$$

	Using \eqref{hyp:pqr} we find that the only two possible points of local extrema of $h(t)$ are 
	$$t_1 = 0 \quad \mbox{and} \quad t_2 = \left( \frac{2^{p+1}(1-qr)(p+1)}{(q+1)(p-r)}\right)^{\frac{q}{pq - 1}}.$$

	Again, inequalities in \eqref{hyp:pqr} yield $h(0) = 0$, $\limit{\toinfty{t}} h(t) = -\infty$ and also $h'$ is positive for $t \ll 1$. Thus, $t_2$ is a point of global maximum of $h$ and $h(t_2) > 0$. We denote $T := t_2$.

	\mezeraA

	Let $v \in X$ satisfy $\normX{v} = T$. Inequality (\ref{J(v) estimated using advanced Poincare}) for $v$ now reads as
	
	\begin{equation}
		J(v) \geq T^{r+1} \left( h(T) - \frac{\lambda}{2^{r+1} (r+1)} \right).
		\label{J(v) estimated using advanced Poincare with T}
	\end{equation}

	Denote
	\begin{equation}
	\label{lambda0}
	\lambda_0 := 2^{r+1}(r+1) \left( \frac{q}{q+1} \, T^{\frac{1}{q} - r} - \frac{1}{2^{p+1}(p+1)} \, T^{p-r} \right) = 2^{r+1}(r+1) \, h(T)
	\end{equation}

	and let $\lambda \in (0,\lambda_0)$. Setting
	$$C := T^{r+1} \left( \frac{q}{q+1} \, T^{\frac{1}{q}-r} - \frac{1}{2^{p+1}(p+1)} \, T^{p-r} - \frac{\lambda}{2^{r+1} (r+1)} \right) = T^{r+1} \left( h(T) - \frac{\lambda}{\lambda_0} \, h(T) \right),$$
	
	it follows from (\ref{J(v) estimated using advanced Poincare with T}) that for $v\in X$ with $\|v\|_X=T$, 
	$$J(v) \geq C > 0,$$
	which proves the claim.
\end{proof}

Let $B_T \subset X$ denote the closed ball centered at the origin with radius $T$, that is,
	$$B_T = \left\{ u \in X: \hspace{1 ex} \normX{u} \leq T \right\}.$$

\begin{lemma} \label{lemma: J(v) negative outside B_T}
	There exists $e \in X\setminus B_T$ such that $J(e) < 0$.
\end{lemma}
\begin{proof}
	Let $\varphi \in X$ be an arbitrary, but fixed function such that $\varphi \in \partial B_T$ and $\varphi > 0$ in $(0,1)$.
	Using \eqref{J(v) norms}, for $t$ positive

\begin{equation}\label{eq: J in t phi}
    J(t \, \varphi) = \frac{q}{q+1} \, t^\frac{q+1}{q} \, \|\varphi\|_X^\frac{q+1}{q} 
	- \frac{\lambda}{r+1} \, t^{r+1} \, \|\varphi\|_{L^{r+1}(0,1)}^{r+1} 
	- \frac{1}{p+1} \, t^{p+1} \, \|\varphi\|^{p+1}_{L^{p+1}(0,1)},
\end{equation}
	and from \eqref{hyp:pqr} and the assumptions for $\lambda$, we have

	$$\limit{\toinfty{t}} J(t \, \varphi) = -\infty.$$

	Therefore for a sufficiently large $t > 1$, the function $e := t \, \varphi$ is such that $e \in X\setminus B_T$ and $J(e) < 0$.
\end{proof}

Now we are in a position to prove the existence of a first solution to~\BVP.

\begin{prop} \label{TH: MP solution}
Let $\lambda_0 > 0$ be as in Lemma \ref{lemma: J(v) positive on a sphere}. 
For $\lambda \in (0, \lambda_0)$, there exists a nontrivial weak solution $v_1 \in X$ of \BVP. Moreover, $J(v_1) > 0$.
\end{prop}
\begin{proof}
Observe that $J \in C^1(X, \R)$ and satisfies the PS-condition (see Section \ref{sect:PS}). Also, $J(0) = 0$ and from Lemma \ref{lemma: J(v) positive on a sphere}, for $\lambda \in (0, \lambda_0)$, $\infimum{v \in \partial B_T} J(v) \geq C > 0$. Finally, from Lemma \ref{lemma: J(v) negative outside B_T}, there exists $e \in X$ such that $\normX{e} > T$ and $J(e) < 0$. Let 
$$
c:=\inf \limits_{\gamma\in \Gamma} \max \limits_{0\leq s\leq 1}J(\gamma(s)),
$$
where $\Gamma:=\{\gamma\in C([0,1];X):\,\gamma(0)=0 \mbox{ and }\gamma(1)=e\}$. The Mountain Pass Theorem (see Theorem 6.4.5 in \cite{drabekmilota} with $F = J$ and $R = T$) yields the existence of $v_1 \in X$ such that 
	\begin{equation}
		\label{MP zaver I.}
		J(v_1) = c	\quad \hbox{and} \quad DJ(v_1)= 0 \quad \hbox{in } X^*.
  \end{equation}
The definition of $c$ yields that $J(v_1) \geq C > 0$. From (\ref{MP zaver I.}), $v_1$ is a nontrivial weak solution of \BVP. This completes the proof.
\end{proof}

\section{Local minimum solution} \label{sec:second_sln}
We recall again that $p,q,r$ satisfy \eqref{hyp:pqr}. In this part we show that for $\lambda>0$ small enough, there exists a second nontrivial solution of \BVP, which is a local minimum for $J$ near the trivial solution. 

We use the notations and conventions introduced in Section \ref{sec:first_sln}. Let $T$, $\lambda_0$ and $C$ be as in the statement of Lemma \ref{lemma: J(v) positive on a sphere} and recall that $B_T$ is the closed ball in $X$ centered at zero with radius $T$.

\begin{lemma} \label{lemma: J is negative in B_T}
	Let $\lambda\in (0,\lambda_0)$. There exists $\tilde{v} \in X\setminus \{0\}$ such that $\tilde{v} \in \interior B_T$ and $J(\tilde{v}) < 0$.
\end{lemma}
\begin{proof}
	Let $\varphi \in X$ be an arbitrary, but fixed function such that $\varphi \in \partial B_T$ and $\varphi > 0$ in $(0,1)$.
	
	Then, using \eqref{eq: J in t phi}, for $t \in (0,1)$ and $\lambda > 0$
	\begin{equation*}
		J(t \, \varphi) = t^{r+1} \, \left( \frac{q}{q+1} \, t^{\frac{1}{q}-r} \, \|\varphi\|_X^\frac{q+1}{q} - \frac{\lambda}{r+1} \, \|\varphi\|_{L^{r+1}(0,1)}^{r+1} - \frac{1}{p+1} \, t^{p-r} \, \|\varphi\|^{p+1}_{L^{p+1}(0,1)} \right). 
	\end{equation*}
	
	Since $\frac{1}{q} > r > 0$, $\lambda \in (0, \lambda_0)$ and since $T>0$ does not depend on $\lambda$, we may select $t=O(\lambda^{\frac{q}{1-qr}})$ small enough such that $t\in (0,T)$ and $J(t \, \varphi) < 0$.
	Finally, if we denote $\tilde{v} := t \, \varphi$, then $\tilde{v} \in B_T$ and $J(\tilde{v}) < 0$.
	This completes the proof of the lemma.
\end{proof}

The existence of another solution is next stated.

\begin{prop} \label{TH: solution near the origin}
	Let $T>0$ and $\lambda_0$ be as in Lemma \ref{lemma: J(v) positive on a sphere}. 
	For all $\lambda \in (0, \lambda_0)$, there exists a nontrivial weak solution $v_2 \in \interior B_T$ of \BVP. In addition, $J(v_2) < 0$.
\end{prop}
\begin{proof}
Consider the minimization problem 
\begin{equation}\label{eqn:min_prob}
   m_\lambda:= \inf \limits_{v\in B_T}J(v).
\end{equation}
We claim that $m_\lambda$ is attained. To prove this we proceed as follows. Recall that $X$ is reflexive\footnote{See Lemma \ref{lemma: norm equivalence, Banach, Reflexive}.} and thus, using Kakutani's Theorem (see \cite{brezis}, Theorem 3.17, p. 67),	$B_T$ is sequentially weakly compact. It suffices to prove that $J$ is sequentially weakly lower semicontinuous in $B_T$. Once this is proven, the existence of $v_2 \in B_T$ such that $J(v_2) =m_{\lambda}$
	will follow from Theorem 1.1 in \cite{struwe}, and the corresponding comments. In particular, $m_{\lambda}\in \R$. Also, Lemma \ref{lemma: J(v) positive on a sphere} and Lemma \ref{lemma: J is negative in B_T} guarantee that $v_2 \in \interior B_T$ with $J(v_2) < 0$. Consequently, $v_2$ is a nontrivial weak solution of \BVP.\\
 
 Now we prove the sequential weak lower semicontinuity of $J$ in $B_T$. Let $(v_n) \subset B_T$ be arbitrary and such that $v_n \rightharpoonup v$ weakly in $X$. Since $B_T$ is sequentially weakly compact, $v \in B_T$. We prove that
	$$J(v) \leq \limesinferior{n \rightarrow +\infty} J(v_n).$$

Notice that $X$ is continuously embedded in $W^{1,\frac{q+1}{q}}(0,1)$ and hence compactly embedded into $C[0,1]$. Thus, the sequence $(v_n)$ converges uniformly to $v$ in $[0,1]$. In particular, $(v_n)_+\to v_+$ uniformly in $[0,1]$. Since the norm $\normX{.}$ is sequentially weakly lower semicontinuous\footnote{From \cite{brezis}, Prop. 3.5 (iii), p. 58} 
	in $B_T$, we estimate
	\begin{equation}
 \begin{aligned}
		\limesinferior{\ninfty} J(v_n) &\geq \frac{q}{q + 1} \, \normX{v}^{\frac{q + 1}{q}}
		- \frac{\lambda}{r+1} \, \lim \limits_{\ninfty} \left( \int\limits_0^1 (v_n)_+^{r+1}  \, \mbox{d}x \right)
		- \frac{1}{p+1} \, \lim \limits_{\ninfty} \left( \int\limits_0^1 (v_n)_+^{p+1} \, \mbox{d}x \right)\\
  &= \frac{q}{q + 1} \, \normX{v}^{\frac{q + 1}{q}}
		- \frac{\lambda}{r+1}\left( \int\limits_0^1 v_+^{r+1}  \, \mbox{d}x \right)
		- \frac{1}{p+1} \left( \int\limits_0^1 v_+^{p+1} \, \mbox{d}x \right)\\
  &= J(v).
		\label{J is almost s.w.l.semicontinuous}
  \end{aligned}
	\end{equation}

This proves the last claim and completes the proof of the proposition. 	
\end{proof}

\section{Palais-Smale condition}
\label{sect:PS}
\begin{lemma} \label{lemma: PSc verified}
The energy functional $J$ satisfies {\rm PS}-condition.
\end{lemma}
\begin{proof}
	Let us consider a sequence $(v_n) \subset X$ satisfying (\ref{Assumption Palais-Smale}) with $c \in \R$.
	This implies that there exists $K > 0$ such that
	\begin{equation}
		J(v_n) \leq K \quad \forall n \in \N.
		\label{P-S -> J is bounded}
	\end{equation}
	
It also implies that for any $\varepsilon > 0$ there exists $n_0 \in \N$ such that for all $n \in \N$ with $n\geq n_0$ and for all $\varphi \in X$ we have
	\begin{equation}
		|DJ(v_n)\varphi| \leq \varepsilon \, \normX{\varphi}.
		\label{P-S -> |DJ| arbitrary small}
	\end{equation}

	\mezeraA

	First we prove that $(v_n)$ is bounded in $X$.
	Passing to a subsequence of $(v_n)$ (which for simplicity we denote the same), arguing by contradiction, we assume $\normX{v_n} \rightarrow +\infty$.
	
	Using (\ref{P-S -> J is bounded}) and (\ref{P-S -> |DJ| arbitrary small}), we estimate
	$$	K + \frac{\varepsilon}{p + 1} \, \normX{v_n} \geq J(v_n) - \frac{1}{p + 1} \, DJ(v_n)v_n.$$
	
	Assumption $pq>1$ and Remark \ref{lemma: L infinity norm <= X norm of u} (using, for simplicity, a weaker estimate $\|i\|\leq 1$) yield
	\begin{equation}
	K + \frac{\varepsilon}{p + 1} \, \normX{v_n} \geq \left(\frac{q}{q + 1} - \frac{1}{p + 1}\right) \, \normX{v_n}^{\frac{q + 1}{q}}
	+ \left( \frac{\lambda}{p + 1} - \frac{\lambda}{r + 1}\right) \, \normX{v_n}^{r + 1}.
	\label{P-S clever subtraction}
	\end{equation}
	
	Since $(v_n)$ is assumed to be unbounded, for all $n$ sufficiently large, $\normX{v_n} > 0$ and therefore we can divide both sides of 
	(\ref{P-S clever subtraction}) by $\normX{v_n}^{\frac{q + 1}{q}}$, which gives
	\begin{equation}
		\frac{K}{\normX{v_n}^{\frac{q + 1}{q}}} + \frac{\varepsilon}{(p + 1) \, \normX{v_n}^\frac{1}{q}} \geq \left(\frac{q}{q + 1} - \frac{1}{p + 1}\right)
		+ \left( \frac{\lambda}{p + 1} - \frac{\lambda}{r + 1}\right) \, \frac{1}{\normX{v_n}^{\frac{1}{q} - r}}.
		\label{P-S clever division}
	\end{equation}
	
	From \eqref{hyp:pqr}, we know that all the powers of $\normX{v_n}$ in (\ref{P-S clever division}) are positive and
	$$\frac{q}{q + 1} - \frac{1}{p + 1} > 0.$$

	Thus, taking the limit as $\ninfty$ in (\ref{P-S clever division}) yields a contradiction and hence $(v_n)$ is bounded.

	Passing to a subsequence, denoted for simplicity as $(v_n)$, from Eberlain-Smulyan's Theorem (see \cite{drabekmilota}, Th. 2.1.25, p. 67) and from Lemma \ref{lemma: norm equivalence, Banach, Reflexive}, and since $X$ is compactly embedded into $C[0,1]$, there exists $v \in X$ such that
	\begin{eqnarray}
		v_n &\rightharpoonup& v \quad \mbox{in } X,\label{P-S -> weak convergence}\\
		v_n &\rightarrow& v \quad \mbox{in } L^{\frac{q+1}{q}}(0,1)\label{P-S -> strong convergence}.
	\end{eqnarray}

	We now prove that $(v_n)$ converges strongly in $X$.
	Since (\ref{P-S -> weak convergence}) holds and $X$ is a uniformly convex space, 
	it follows from \cite{brezis}, Prop. 3.32, p. 78, that it suffices to show that $\normX{v_n} \rightarrow \normX{v}$.

	Denote 
	$$\varepsilon_n := \|DJ(v_n)\|_{X^*} := \supremum{\|\varphi\|_X = 1} |DJ(v_n) \, \varphi|.$$
	Observe that $\varepsilon_n \geq 0$ and $\limit{\ninfty} \varepsilon_n = 0$.	
	Choosing $\varphi := v_n - v$ for any $n \in \N$, (\ref{P-S -> |DJ| arbitrary small}) reads as
	\begin{equation*}
		\begin{aligned}
			|DJ(v_n)\varphi| = \int\limits_0^1 |v_n''|^{\frac{1}{q}-1} \, v_n'' \, (v_n'' - v'') \mbox{d}x - \int\limits_0^1 \left( \lambda \, (v_n)_+^{r-1} \, v_n \, 
			(v_n - v) + (v_n)_+^{p-1} \, v_n \, (v_n - v) \right) \mbox{d}x \leq\\
			\leq \varepsilon_n \, \normX{\varphi}.\\
		\end{aligned}
	\end{equation*}
	
	Using the triangle inequality, we estimate
	\begin{equation} \label{P-S -> dosazeno v_n - v}
		\begin{aligned}
			\left| \int\limits_0^1 |v_n''|^{\frac{1}{q}-1} \, v_n'' \, (v_n'' - v'') \, \mbox{d}x\right| - \left| \int\limits_0^1 \lambda \, (v_n)_+^{r-1} \, v_n \, (v_n - v) \, \mbox{d}x \right| -\\
			\qquad - \left| \int\limits_0^1 (v_n)_+^{p-1} \, v_n \, (v_n - v) \, \mbox{d}x \right| \leq \varepsilon_n \, \normX{v_n - v}.
		\end{aligned}
	\end{equation}
	
	Using that $(v_n)_+^{r}\leq |v_n|^{r}$ and the Hölder inequality for the conjugate exponents $\frac{q+1}{q}$ and $q+1$ yields
	\begin{equation*}
	\left| \int\limits_0^1 \lambda \, (v_n)_+^{r-1} \, v_n \, (v_n - v) \, \mbox{d}x \right| \leq |\lambda| \, \normLp{v_n}{r(q+1)}^r \, \normLp{v_n - v}{\frac{q+1}{q}}.
	\end{equation*}
	
	Since $(v_n)$ is a bounded sequence in $X$, from Remark \ref{lemma: L infinity norm <= X norm of u}, and (\ref{P-S -> strong convergence}), we obtain
	\begin{equation}
		\limit{\ninfty} \int\limits_0^1 \lambda \, (v_n)_+^{r-1} \, v_n \, (v_n - v) \, \mbox{d}x = 0.
		\label{P-S -> lim1}
	\end{equation} 
	
	Arguing in a similar manner yields 
		\begin{equation}
	\limit{\ninfty} \int\limits_0^1 (v_n)_+^{p-1} \, v_n \, (v_n - v) \, \mbox{d}x = 0.
	\label{P-S -> lim2}
	\end{equation} 
	
	We recall that $(v_n)$ is bounded in $X$ and $\varepsilon_n \rightarrow 0$. Therefore, expressions  (\ref{P-S -> dosazeno v_n - v}), (\ref{P-S -> lim1}), and (\ref{P-S -> lim2}) imply
	\begin{equation}
		\limit{\ninfty} \int\limits_0^1 |v_n''|^{\frac{1}{q}-1} \, v_n'' \, (v_n'' - v'') \, \mbox{d}x = 0.
		\label{P-S -> lim3}
	\end{equation}
	
	In addition, from the definition of the weak convergence given by (\ref{P-S -> weak convergence}), we know that 
	\begin{equation}
		\limit{\ninfty} \int\limits_0^1 |v''|^{\frac{1}{q}-1} \, v'' \, (v_n'' - v'') \, \mbox{d}x = 0
		\label{P-S -> lim4}
	\end{equation}
	
	and hence subtracting terms in (\ref{P-S -> lim3}) and (\ref{P-S -> lim4}) and using the Hölder inequality yields
	\begin{equation*}
		\begin{aligned}
			0 &= 
			 \limit{\ninfty} \left( \normX{v_n}^{\frac{q+1}{q}} - \int\limits_0^1 |v_n''|^{\frac{1}{q}-1} \, v_n'' \, v'' \, \mbox{d}x - \int\limits_0^1 |v''|^{\frac{1}{q}-1} \, v'' \, v_n'' \, \mbox{d}x + \normX{v}^{\frac{q+1}{q}}
			\right)\\
			&\geq \limit{\ninfty} \left( \normX{v_n}^{\frac{q+1}{q}} - \normX{v_n}^{\frac{1}{q}} \, \normX{v} - \normX{v}^{\frac{1}{q}} \, \normX{v_n} + \normX{v}^{\frac{q+1}{q}}
			\right)\\
			&= \limit{\ninfty} \left( \normX{v_n}^\frac{1}{q} - \normX{v}^\frac{1}{q} \right) \, \left( \normX{v_n} - \normX{v} \right).
		\end{aligned}
	\end{equation*}
	
	Since the function $x \mapsto x^{\frac{1}{q}}$ is strictly increasing, 
	\begin{equation*}
		0 \geq  \limit{\ninfty} \left( \normX{v_n}^\frac{1}{q} - \normX{v}^\frac{1}{q} \right) \, \left( \normX{v_n} - \normX{v} \right) \geq 0,
	\end{equation*}
	
	thus, necessarily,
	\begin{equation}
		\normX{v_n} \rightarrow \normX{v}.
		\label{P-S -> konvergence norem}
	\end{equation}
	
	Assertions (\ref{P-S -> weak convergence}) and (\ref{P-S -> konvergence norem}) prove the statement.
\end{proof}

\section{Proof of Theorem \ref{theo:main_theo_1} and range of $\lambda$}\label{sec:proof_of_main_theorem}
We continue the theoretical discussion with the proof of our main result.

\begin{proof}[Proof of Theorem \ref{theo:main_theo_1}]
	Let $\lambda_0$ be as in Lemma \ref{lemma: J(v) positive on a sphere} and let us consider $\lambda \in (0, \lambda_0)$.
	Then it follows from Propositions \ref{TH: MP solution} and \ref{TH: solution near the origin} that there exist two weak nontrivial solutions 
	$v_1, v_2 \in X$ of \BVP. Moreover, since $J(v_1) > 0 > J(v_2)$, the weak solutions are necessarily distinct. Using Proposition \ref{theo:regl_4th_order_BVP}, we verify that $v_1, v_2$ are classical solutions of \BVP.  Expression (\ref{u vyjadreno pomoci v}) yields the corresponding nontrivial smooth functions $u_1, u_2$ such that the pairs of functions $(u_1, v_1)$ and
	$(u_2, v_2)$ solve (\ref{BVP positive parts}). As it was described 
	at the beginning of Section \ref{sec:setting},
	all these functions are necessarily
	nonnegative. In conclusion, for $\lambda \in (0, \lambda_0)$, the pairs $(u_1, v_1)$ and $(u_2, v_2)$ represent two distinct nontrivial nonnegative classical solutions to \eqref{BVP}-\eqref{eqn:BVP_Dirch}. This concludes the proof of the theorem. 
\end{proof}

Theorem \ref{theo:main_theo_1} states the existence of at least two distinct solutions of (\ref{BVP}) for $\lambda \in (0, \lambda_0)$. This means that in the bifurcation diagram, we find at least two branches of solutions for small positive values of $\lambda$. A natural question concerns with the maximum value of $\lambda_0$ for which Theorem \ref{theo:main_theo_1} is still valid.\\ 

As we will see in Section \ref{sec:numerical_exp}, when $p = 3$, $q = 1.5$ and $r = 3^{-1}$, the numerical illustration predicts the existence of the upper and lower branches for
$\lambda \in (0, \lambda_{\rm bif})$, where $\lambda_{\rm bif} \approx 49$.
For the same values of the parameters $p$, $q$, and $r$, relation \eqref{lambda0} provides us
\begin{equation*}
	\lambda_0 \approx 2.21 \ll \lambda_{\rm bif}.
\end{equation*}

This means that the theoretical results herein describe the system only in a narrow interval for $\lambda$. In part, this non-optimality for the range of $\lambda$ is due to the non-optimality of the constants in the embeddings stated in Lemma \ref{lemma: iterated Morrey's} and Remark \ref{lemma: L infinity norm <= X norm of u}. 

In order to extend the range of values of $\lambda$ for which Theorem \ref{theo:main_theo_1} is still valid, we track up the energy estimates in the proof of Lemma \ref{lemma: J(v) positive on a sphere}. At this point it is reasonable to make use of the optimal constant for the embedding $X \hookrightarrow L^s(0,1)$ in Remark \ref{lemma: L infinity norm <= X norm of u}. For instance, similarly to the well-known result for the Laplace operator, it is readily verified that the infimum 
\begin{equation}
\label{eqn:C_emb}
	C_{emb}^{-1} := \infimum{\substack{u \in X\\u \neq 0}} \frac{\normX{u}}{\normLp{u}{\frac{q+1}{q}}}
\end{equation}
is a positive number and that it is attained. It can also be shown that $C_{emb}^{-1}$ corresponds to the principal eigenvalue of the problem
\begin{equation*}
	\left\{ 
	\begin{aligned}
		\frac{\mbox{d}^{2}{\ }}{\mbox{d}{x^2}} \left( |u''|^{\gamma - 2} \, u'' \right) &= \lambda \, |u|^{\gamma - 2} \, u,\quad x \in (0,1),\\ 
		u(0) = u(1) &= u''(0) = u''(1) = 0.
	\end{aligned}
	\right.
\end{equation*}

We do not know the exact value of this principal eigenvalue, however, in \cite{benedikt} it is proved that
\begin{equation}\label{EQ: estimate benedikt}
	C_{emb} \leq K_{emb} := \frac{\left(\frac{1}{2}\right)^2}{2} \, \min \left\{ \left( \frac{\sqrt{\pi} \, \Gamma 
	(\gamma)}{\Gamma (\gamma + \frac{1}{2})} - \frac{1}{\gamma} \right)^{\frac{1}{\gamma}}, 
	\left( \frac{\sqrt{\pi} \, \Gamma (\gamma')}{\Gamma (\gamma' + \frac{1}{2})} - \frac{1}{\gamma'} 
	\right)^{\frac{1}{\gamma'}}\right\},
\end{equation}
where $\gamma=\frac{q+1}{q}$, $\gamma' := \frac{\gamma}{\gamma - 1}$ and $\Gamma(z) := \int\limits_0^{+\infty} t^{z-1} \, e^{-t} \, \mbox{d}t$.

Let us now improve the estimate on $\lambda_0$. Since $qr<1$, for any $v \in X$
\begin{equation*}
	\normLp{v_+}{r+1} \leq \normLp{v}{r+1} \leq \normLp{v}{\frac{q+1}{q}}
\end{equation*}
and using \eqref{eqn:C_emb} and \eqref{EQ: estimate benedikt}, we obtain that
\begin{equation}\label{EQ: better estimate for Lr norm}
	\normLp{v_+}{r+1} \leq K_{emb} \, \normX{v}.
\end{equation}

From (\ref{J(v) norms}) we estimate
\begin{equation}\label{EQ: J estimated from below better emb}
\begin{aligned}
     J(v) &\geq \frac{q}{q+1} \, \normX{v}^{\frac{q+1}{q}} - \frac{\lambda \, K_{emb}^{r+1}}{r+1} \, \normX{v}^{r+1} 
	- \frac{1}{2^{p+1} (p+1)} \, \normX{v}^{p+1}\\
&\geq \normX{v}^{r+1} \left( \frac{q}{q+1} \, \normX{v}^{\frac{1}{q}-r} - \frac{\lambda \, K_{emb}^{r+1}}{r+1} 
	- \frac{1}{2^{p+1} (p+1)} \, \normX{v}^{p-r} \right).
 \end{aligned}
\end{equation}

If we replace (\ref{J(v) estimated using advanced Poincare}) by (\ref{EQ: J estimated from below better emb}) in the proof
of Lemma \ref{lemma: J(v) positive on a sphere}, we get\footnote{For definition of $T$, see the proof 
of Lemma \ref{lemma: J(v) positive on a sphere}.}
\begin{equation}\label{EQ: better lambda0}
	\lambda_0 = \frac{r+1}{K_{emb}^{r+1}} \left( \frac{q}{q+1} \, T^{\frac{1}{q}-r} - \frac{1}{2^{p+1} (p+1)} \, T^{p-r} \right).
\end{equation}

When $p = 3$, $q = 1.5$, and $r = 3^{-1}$, we approximate this updated value of $\lambda_0$ as $\lambda_0 \approx 16.02$.
Even though, the previous procedure improves the value of $\lambda_0$, the estimate is still far from optimal. One of the reasons is that the optimality was used only for the embedding $X \hookrightarrow L^{\frac{q+1}{q}}(0,1)$ to treat the term with $L^{r+1}$-norm. The term with $L^{p+1}$-norm was treated as before. Another reason is that solutions predicted in Proposition \ref{TH: MP solution} have positive energy, but numerical simulations in Section \ref{sec:numerical_exp} indicate there exist solutions of Mountain pass type with nonpositive energy.\\

\section{Symmetry of solutions} \label{section:Symmetry}

Our next result states symmetry properties of solutions of \eqref{BVP}-\eqref{eqn:BVP_Dirch}.

\begin{prop}\label{theo:symmetry_slns}
Let $\lambda > 0$ and let $(u,v)$ be a classical solution of the system \eqref{BVP}-\eqref{eqn:BVP_Dirch} with $u,v>0$ in $(0,1)$. Then,
\begin{itemize}
\item[i.] $u$ and $v$ are symmetric with respect to the vertical line $x=\frac{1}{2}$;
\item[ii.] $u\big(\frac{1}{2}\big)=\|u\|_{L^{\infty}(0,1)}$ and $v\big(\frac{1}{2}\big)=\|v\|_{L^{\infty}(0,1)}$;
\item[iii.] $x=\frac{1}{2}$ is the only critical point of $u$ and $v$.
\end{itemize}
\end{prop}

The proof of Proposition \ref{theo:symmetry_slns} follows the (by now) standard method of \emph{moving planes} (see  \cite{berestycki1988monotonicity, gidas1979symmetry} and Section 9 in \cite{de2008semilinear}). Nonetheless, we present a self-contained proof, adapted to the one dimensional case treated in this work.\\ 

We remark in advance that we only require the monotonicity of the nonlinearities and hence in this part it is enough to assume  $0<r,p,q<\infty$ and  $\lambda$ to be nonnegative. Since we deal with the one dimensional case, the proof proceeds with an ad-hoc procedure described by a series of lemmas. The core idea is based upon the classic method of moving planes (see \cite{de1994monotonicity},\cite{gidas1979symmetry}).
  
\vskip 3pt

\begin{lemma}\label{lemma:concv_critc_pt}
Let $-\infty\leq a<b\leq +\infty$ and $w\in C^2(a,b)$ be such that $w''$ does not vanish in $(a,b)$, then $w$ has at most one critical point in $(a,b)$.
\end{lemma}

\begin{proof}
It follows directly by a direct application of the \emph{Mean Value Theorem} to $w'$.
\end{proof}

\vskip 6pt
\begin{lemma}\label{lemma:hopf's_lemma} Let $-\infty<a<b<+\infty$ and $w\in C^2(a,b)\cap C^1[a,b]$. Then 
\begin{itemize}
\item[i.] $w''\leq 0$ in $(a,b)$, $w(a)=0$ and $w(b)\geq 0$ implies that either $w\equiv 0$ in $[a,b]$ or
$$
w>0 \quad \hbox{in}\quad (a,b) \quad \hbox{and}\quad  w'(a)>0.
$$

\item[ii.]$w''\geq 0$ in $(a,b)$, $w(a)=0$ and $w(b)\leq 0$ implies that either $w\equiv 0$ in $[a,b]$ or 
$$
w<0 \quad \hbox{in}\quad (a,b) \quad \hbox{and}\quad  w'(a)<0.
$$
\end{itemize}
\end{lemma}

\begin{proof}
We only prove i., since ii. is analogous. 
First, let $x_0\in (a,b]$ be arbitrary, but fixed. Consider $z(x):=w(x) - \frac{w(x_0)}{x_0-a}(x-a)$ for $x\in [a,x_0]$. Observe that
$z''=w''\leq 0$ in $(a,x_0)$, $z(a)=w(a)=0$ and $z(x_0)=0$. The concavity of $z$ implies that $z\geq 0$ in $(a,x_0)$. Thus, $w(x)\geq  \frac{w(x_0)}{x_0-a}(x-a)$ in $[a,x_0]$. In particular, from the definition of derivative, $w'(a)\geq\frac{w(x_0)}{x_0-a}$. 

\vskip 3pt
Proceeding similarly one can show that $w(x)\geq w(b)-\frac{w(b)-w(x_0)}{b-x_0}(b-x)\geq 0$ in  $[x_0,b]$ and $w'(b)\leq \frac{w(b)-w(x_0)}{b-x_0}\leq 0$. Taking $x_0=b$, we conclude that $w\geq 0$ in $[a,b]$.

\vskip 3pt
Now, let $x_0\in [a,b]$ such that $w(x_0)=\|w\|_{L^{\infty}(a,b)}$. If $w(x_0)=0$, then $w\equiv 0$ in $[a,b]$. Otherwise,  $x_0\in (a,b]$ and $w(x_0)>0$. The above developments imply that $w(x)\geq \frac{w(x_0)}{x_0-a}(x-a)>0$ in $(a,x_0)$, $w'(a)>\frac{w(x_0)}{x_0-a}$ and $w(x)\geq w(b)-\frac{w(b)-w(x_0)}{b-x_0}(b-x)>0$ in  $[x_0,b]$. This completes the proof of the lemma. 
\end{proof}

\begin{rmk}\label{rmk:Hopf'slemma}
From the proof of Lemma \ref{lemma:hopf's_lemma}, in the case that $w'' \leq 0$ in $(a,b)$, $w(a)=w(b)=0$ and $\|w\|_{L^{\infty}(a,b)}>0$, then $w'(a)>0>w'(b)$.
\end{rmk}

Next, consider the space
\begin{equation}\label{def:X_0}
X_0:=C^2(0,1)\cap \big\{ w\in C^1[0,1]\,:\,w(0)=w(1)=0 \big\}.
\end{equation}

For $w\in X_0$ with $w'(0)>0$ and $w'(1)<0$, set
$$
\begin{aligned}
A_w:=&\sup \big\{a\in (0,1)\,:\, w'>0 \quad \hbox{in} \quad (0,a)\big\}\\
B_w:=&\inf \big\{b\in (0,1)\,:\, w'<0 \quad \hbox{in} \quad (b,1)\big\}.
\end{aligned}
$$

Observe that $A_w,B_w$ are well defined and $A_w\in (0,1]$ and $B_w\in [0,1)$. Since $w\in X_0$, $A_w, B_w\in (0,1)$.

\vskip 6pt
\begin{lemma}\label{lemma:AwBw_crit_points}
Let $w\in X_0$ be as above. Then $A_w$ and $B_w$ are critical points of $w$. 
\end{lemma}

\begin{proof}
 By the approximation property of the supremum and the infimum, we find that $w'(A_w)\geq 0$ and $w'(B_w)\leq 0$. If $w'(A_w)>0$, the continuity of $w'$ yields that for some $\delta>0$ small, $w'>0$ in an interval of the form $(0,A_w+\delta)$, thus violating the definition of $A_w$. We conclude that $w'(A_w)=0$. Similarly, we prove that $w'(B_w)=0$.   
\end{proof}

\vskip 6pt
\begin{proof}[Proof of Proposition \ref{theo:symmetry_slns}] Let $u,v\in C^2(0,1)\cap C^1[0,1]$ with $u,v>0$ in $(0,1)$ and such that the pair $(u,v)$ is a classical solution of \eqref{BVP}-\eqref{eqn:BVP_Dirch}. Let $\delta \in (\frac{1}{2},1)$ be arbitrary, but fixed. Observe that $0<2\delta -1<1$. Write
$$
{u}_{\delta}(x):=u(2\delta-x) \quad \hbox{and} \quad {v}_{\delta}(x):=v(2\delta -x) \quad \hbox{for} \quad x\in [2\delta-1,1].
$$

Notice that ${u}_{\delta}$ and ${v}_{\delta}$ are the reflections of $u$ and $v$ respectively, with respect to the vertical line $x=\delta$. Also, $u_{\delta}$ and $v_{\delta}$ solve the system \eqref{BVP} in $(2\delta-1,1)$.

\vskip 6pt
Now define
$$
w_{\delta}(x):={u}_{\delta}(x)-u(x) \quad \hbox{and} \quad z_{\delta}(x):={v}_{\delta}(x)-v(x) \quad \hbox {for }x\in [2\delta-1,1].
$$

In view of Lemma \ref{lemma:hopf's_lemma} and Remark \ref{rmk:Hopf'slemma}, for $\delta\in (\frac{1}{2},1)$ close enough to $1$, $u,v$ are strictly decreasing in $(2\delta-1,1)$. Consequently, $w_{\delta},z_{\delta}>0$ in $(\delta,1]$. 

\vskip 3pt
Set
$$
\delta_*:=\inf\Big\{\delta\in\Big(\frac{1}{2},1\Big)\,:\,z_{\delta}>0 \quad \hbox{in} \quad (\delta,1)\Big\}.
$$

We claim that
$$
\delta_*=\inf\big\{\delta\in\Big(\frac{1}{2},1\Big)\,:\,w_{\delta}>0 \quad \hbox{in} \quad (\delta,1)\big\}.
$$

To prove the claim, let $\delta \in (\frac{1}{2},1)$ be such that $z_{\delta}\geq 0$ in $[\delta,1]$. Using \eqref{BVP},
$$
-w_{\delta}''=\lambda(v_{\delta}^r -v^r) +v^p_{\delta}-v^p \geq 0 \quad \hbox{in} \quad (\delta,1).
$$

Since $w_{\delta}(\delta)=0$ and $w_{\delta}(1)=u(2\delta-1)>0$,  Lemma \ref{lemma:hopf's_lemma} implies that $w_{\delta}>0$ in $(\delta,1)$. This proves that $w_{\delta}>0$ in $(\delta,1)$, whenever $z_{\delta}\geq 0$ in $(\delta,1)$. 

\vskip 3pt
Proceeding similarly, $z_{\delta}>0$ in $(\delta,1)$, whenever $w_{\delta}\geq 0$ in $(\delta,1)$. Since $\delta \in (\frac{1}{2},1)$ is arbitrary, the previous discussion proves the claim.

\vskip 3pt
Now, from the definition of $\delta_*$, 
$$
w_{\delta_*}(\delta_*)=z_{\delta_*}(\delta_*)=0.
$$ 

We prove next that $\delta_*=\frac{1}{2}$. Assume by contradiction that $\delta_*>\frac{1}{2}$ and notice that $w_{\delta_*},z_{\delta_*}\geq 0$ in $(\delta_*,1)$. 

\vskip 3pt
Also, since $w_{\delta_*}(1)=u(2{\delta_*}-1)>0$ and $z_{\delta_*}(1)=v(2{\delta_*}-1)>0$, Lemma \ref{lemma:hopf's_lemma} yields that $w_{\delta_*},z_{\delta_*}>0$ in $(\delta_*,1)$. 

The continuity of the family of functions $\{w_{\delta}\}_{\delta}$ and $\{z_{\delta}\}_{\delta}$ in $C^1[0,1]$ with respect to the parameter $\delta$, allows us to find $\hat{\delta}\in (\frac{1}{2},\delta_*)$ such that $w_{\hat{\delta}},z_{\hat{\delta}}>0$ in $(\hat{\delta},1)$. This contradicts the definition of $\delta_*$ and proves that $\delta_*=\frac{1}{2}$.

\vskip 3pt
Since $w_{\frac{1}{2}}(x)=u(1-x)-u(x)\geq 0$ and $z_{\frac{1}{2}}(x)=v(1-x)-v(x)\geq 0$ for every $x\in [\frac{1}{2},1]$, we find that $u(1-x)\geq u(x)$, $v(1-x)\geq v(x)$ for any $x\in [\frac{1}{2},1]$. A similar argument shows that $u(1-x)\geq u(x)$, $v(1-x)\geq v(x)$ for any $x\in [0,\frac{1}{2}]$ and consequently for $x\in[0,1]$.

\vskip 3pt
Now notice that $u(1-x),v(1-x)$ also solve \eqref{BVP}-\eqref{eqn:BVP_Dirch}. Since $(u_{\frac{1}{2}})_{\frac{1}{2}}=u$ and $(v_\frac{1}{2})_{\frac{1}{2}}=v$, we may argue as above to find that 
$u(1-x)\leq u(x)$, $v(1-x)\leq v(x)$ for any $x\in [0,1]$. Thus, $u,v$ are symmetric with respect to $x=\frac{1}{2}$.

\vskip 3pt
Now, let $x_u,x_v\in [0,1]$ be such that 
$$
u(x_u)=\|u\|_{L^{\infty}(0,1)} \quad \hbox{and} \quad  v(x_v)=\|v\|_{L^{\infty}(0,1)}.
$$ 

Since $u(x_u)=\max\{u(x)\,:\,x\in [0,1]\}$ and $u'(0)>0>u'(1)\neq 0$, $x_u\in (0,1)$ and consequently $u'(x_u)=0$. Similarly, $x_v\in (0,1)$ and $v'(x_v)=0$.

\vskip 3pt
From Lemmas \ref{lemma:concv_critc_pt} and \ref{lemma:AwBw_crit_points} and the symmetry of $u$ and $v$, $A_u=B_u=x_u=\frac{1}{2}$ and $A_v=B_v=x_v=\frac{1}{2}$ and finishes the proof of the lemma.
\end{proof}

\section{Numerical illustration} 
\label{sec:numerical_exp}

To support the theoretical results in this paper and to obtain wider intuition about behavior of our system, this section discusses the numerical strategy for obtaining visual description of our system.

Let us begin with the strategy for the implementation. First, we transform the system \eqref{BVP} into the associated first order system (of four equations)
\begin{equation}\label{eqn:firstord_ODE}
\left\{ 
		\begin{aligned}
			u'(x) &= w(x),\\
			v'(x) &= z(x),\\
		    w'(x) &= -\lambda \, (v_+(x))^r - (v_+(x))^p,\\
      z'(x) &= -|u(x)|^{q-1} \, u(x),
		\end{aligned}
		\right.
	\end{equation}
and consider this system with the initial conditions	
\begin{equation}\label{eqn:firstord_ODE_IC}
u(0) = 0,\quad v(0) = 0,\quad w(0) = du_0, \quad z(0) = dv_0,    
\end{equation}
where $du_0$ and $dv_0$ are considered as free parameters. Let us assume that the \emph{initial boundary value problem} \eqref{eqn:firstord_ODE}-\eqref{eqn:firstord_ODE_IC} is such that for $(du_0,dv_0)$ in a rectangle $R\subset (0,+\infty)\times (0,+\infty)$, existence, uniqueness and continuity of solutions with respect to the parameters $(du_0,dv_0)$ hold true in $[0,1]$. Then given $(du_0,dv_0)\in R$, the functions $u$, $v$, $w$ and $z$ are continuous in $[0,1]$ and we write
$$
u(x)=u(x,du_0,dv_0) \quad \hbox{and} \quad v(x)=v(x,du_0,dv_0).
$$

We are therefore interested in  the zeroes of the  continuous mapping $\Phi:R\to \R^2$ defined by 
\begin{equation}\label{def:shooting_map_Phi}
\Phi(du_0,dv_0)=(u(1),v(1)). 
\end{equation}

Analyzing the sign changes of the components of $\Phi$, with respect to the values of $du_0$ and $dv_0$, we locate numerically two solutions for a wide range of values $\lambda$.
As mentioned in the Introduction such implementation is motivated by a combination of the standard \emph{Shooting method} and the heuristics of the \emph{Poincar\'{e}-Miranda Theorem}. \\

The strategy is as follows.
\begin{itemize}
    \item Choose $du_0$ and $dv_0$ appropriately, i.e, so that existence and uniqueness in $[0,1]$ and continuity with respect to initial conditions holds.
    
    \item Compute the corresponding solutions $(u,v,w,z)$ of (\ref{eqn:firstord_ODE})-\eqref{eqn:firstord_ODE_IC} in the interval $[0,1]$.
    
    \item Extract the values $u(1)$ and $v(1)$. Observe that for certain values of $du_0$ and  $dv_0$ the corresponding pair of functions $(u,v)$ is a solution to \eqref{BVP positive parts} provided $u(1) = v(1) = 0$.

    \item Given a tolerance $\epsilon>0$,  the numerical calculations of $u$ and $v$ yield an admissible numerical solution to \eqref{BVP positive parts} provided $|u(1)|,|v(1)|<\epsilon$.

    \item It is more efficient to track simultaneously the different regions where either $|u(1)|\geq \epsilon$ or $|v(1)|\geq \epsilon$. This can be interpreted as tracking the sign-changes of the values $u(1)$, $v(1)$ as the parameters $du_0$ and $dv_0$ vary.
\end{itemize}

A reader can imagine the process as if \emph{two football players} kick simultaneously two balls from the ground (zero height) on the left border of the field, each of the kicks performed with a given initial slope. The trajectory of the kicked balls are ruled by the equations in \eqref{BVP positive parts}. The football players aim at a bin located on the right border of the field. To hit the bin with a ball it is necessary to choose the initial slope of the kick so that the ball descends to a zero height exactly on the right border of the field. Since each of the football players may kick with different strength, the required slopes for the football players need not be the same.

The implementation initially uses a ``coarse" grid of different parameters $du_0,dv_0$ to roughly locate the sign-changes. Then, an adaptive strategy is employed in order to improve the grid density. This strategy iterates the values of the parameters $du_0$ and $dv_0$, computes again the corresponding solutions of \eqref{eqn:firstord_ODE} and tracks the regions of the corresponding changes of sign.

We implemented the idea of the experiments as a set of scripts in Matlab. The script is optimized to provide the results in reasonable time and accuracy.

	\begin{description}[itemindent=0pt]
		\item [$<$numerical.m$>$]\ \\
		Define algorithm settings and values of the parameters $p, q, r$.\\
		Define the initial range for $du_0$ and $dv_0$.\\
		Define the range for the parameter $\lambda$.\\
		Iterate through the range for $\lambda$ and do the following:
		
		\begin{description}[itemindent=0pt]
			\item[$<$for cycle$>$]\ \\
			Run \emph{shooting.m} with desired settings.
		
			\begin{description}[itemindent=0pt]
				\item[$<$shooting.m$>$]\ \\
				Represent $du_0$-$dv_0$ plane by coarse ($\Delta = 0.1$) and dense ($\Delta = 0.005$) grids.\footnote{
					The grid is rectangled and uniform, and the distance (in either direction $du_0$ and $dv_0$) between two neighboring vertices is $\Delta$. The results are computed 
					only in the vertices of the grid. For instance, consider the rectangle $[a,b]\times [c,d]\subset (0,\infty)\times (0,\infty)$ for selecting the pairs $(du_0,dv_0)$. The corresponding grid of vertices reads as 
					$$\left\{(a + i \, \Delta, b + j \, \Delta ) \in [a,b]\times[c,d]:\quad i = 0, 1, \ldots, \frac{b-a}{\Delta},\ j = 0, 1, \ldots, \frac{d-c}{\Delta}\right\}.$$ 
					As an example, if $b-a, d-c, \Delta^{-1} \in \N$, then for any unit square in $du_0$-$dv_0$ plane, the coarse grid contains $10^2$ vertices and the dense grid 
					contains $200^2$ vertices.}\\
				Iterate vertices of the coarse grid and run \emph{shootandsolve.m}.
				
				\begin{description}[itemindent=0pt]
					\item[$<$shootandsolve.m$>$]\ \\
					Run \emph{ode45.m} (Runge-Kutta method for ODEs from Matlab library) for \eqref{eqn:firstord_ODE}-\eqref{eqn:firstord_ODE_IC} and any given 
					initial condition $(u(0), w(0), v(0), z(0)) = (0, du_0, 0, dv_0)$.\\
					Compute the solution and calculate the residues $u(1)$ and $v(1)$.\\
					Return the residues as the return value.
					\item[$<$/shootandsolve.m$>$]
				\end{description}
				
				Based on the residues for a pair $(du_0, dv_0)$, assign the pair a color using the following scheme:
                \begin{equation} \label{color legend}
                    \begin{aligned}
        					&u(1) > 0,& &v(1) > 0& &- \qquad \mbox{green};&\\
        					  &u(1) >0 ,& &v(1) < 0& &- \qquad \mbox{yellow};&\\
        					&u(1) < 0,& &v(1) > 0& &- \qquad\mbox{blue};&\\
        					  &u(1) < 0,& &v(1) < 0& &- \qquad \mbox{red}.&
                    \end{aligned}
                \end{equation} 
                
                 Plot the colors in a $du_0$-$dv_0$ diagram and save them into a variable for the coarse grid.
				
				The pair $(du_0, dv_0)$ such that $(u, v)$ is also a solution of (\ref{BVP positive parts}) 
				can be located exactly at the point where all colors meet, in other words, where both residues are zero.
				
				Iterate through the vertices of the coarse grid and choose only the points which have a~neighbour of a different color (we target edges between two colors).
				
				In the dense grid, proceed only with the vertices corresponding to the chosen points in the coarse grid.

				Run \emph{shootandsolve.m} for the vertices in the dense grid.

				Plot the colors in the $du_0$-$dv_0$ diagram and save them into a variable for the dense grid.

				Explore neighbourhood of the points in the dense grid and check whether all colors are present in the neighbourhood
				{(if so, assume there is a solution in the neighbourhood)}.

				Approximate the values of $(du_0,dv_0)$ corresponding to the solution -- denote them by $(solU, solV)$ -- and mark them in the graph with a black circle.

				Return $(solU, solV)$ -- or \emph{(Inf, Inf)} if no solution was found.

				\item[$<$/shooting.m$>$]
			\end{description}

			If the return value is \emph{(Inf, Inf)} (solution not found), exit the {``for''} cycle.

			Run \emph{showsolution.m}, pass $(solU, solV)$ as a parameter.
			\begin{description}[itemindent=0pt]
				\item[$<$showsolution.m$>$]\ \\
				Compute the solution using $(solU, solV)$ and \emph{ode45.m}.

				Plot $\|v\|_{L^{\infty}(0,1)}$ vs $\lambda$, where $v$ is the component $v$ found above, i.e., the numerical solution of \eqref{BVP_merge_eq}.

				Return $L^{\infty}$-norm of the plotted function.

				\item[$<$/showsolution.m$>$]
			\end{description}
			Save the value of $\lambda$ and the corresponding $\|v\|_{L^{\infty}(0,1)}$ of the solution $v$ into a text file.

			Considering the development of the results for the previous values of $\lambda$, automatically adapt the ranges of $du_0$ and $dv_0$ for the next iteration
			(for the optimization, it is necessary to use the narrowest range of the parameters as possible).
			\item[$<$/for cycle$>$]
		\end{description}
	
		Load results for the whole range of $\lambda$ from the text file.

		Vizualize dependancy of the $L^{\infty}$-norm of the solution on the parameter $\lambda$ -- plot the bifurcation diagram.		
		\item [$<$/numerical.m$>$]
	\end{description}

	The use of a coarse grid first and then a denser one in the implementation saves significant amount of time and memory. In the graph, the optimized script skips parts of the domain
	where the results cannot be located (from the coarse grid's point of view), thus the script leaves blank rectangles in the graph.

\vskip 3pt

	Next, we discuss our numerical findings in more detail. Recall that we have set $p = 3$, $q = 1.5$, and $r = \frac{1}{3}$. 
	Nevertheless, for small positive values of $\lambda$, the numerical experiments anticipate existence of a nontrivial solution as it can be seen in Figure \ref{fig: shooting, lambda = 1 LOWER}. 
	\begin{figure}[phtb!]
		\centering
		\begin{subfigure}{0.4\textwidth}
			\centering
			\includegraphics[width = \linewidth]{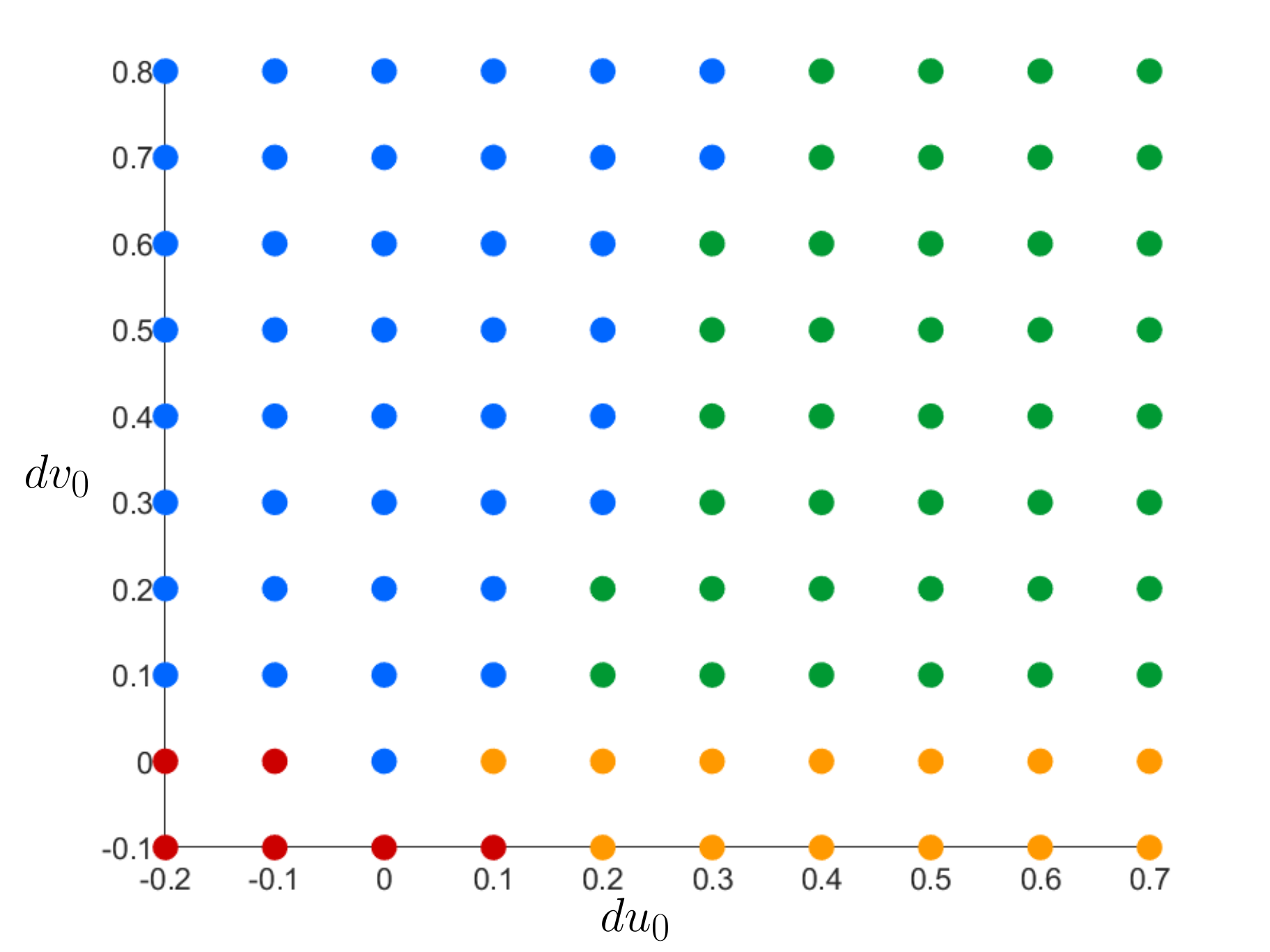}				
			\caption{Diagram for $\lambda = 1$ (coarse grid)\\\ \\ \ }
			\label{fig: shooting, lambda = 1, coarse}
		\end{subfigure} 
		\begin{subfigure}{0.4\textwidth}
			\centering
			\includegraphics[width = \linewidth]{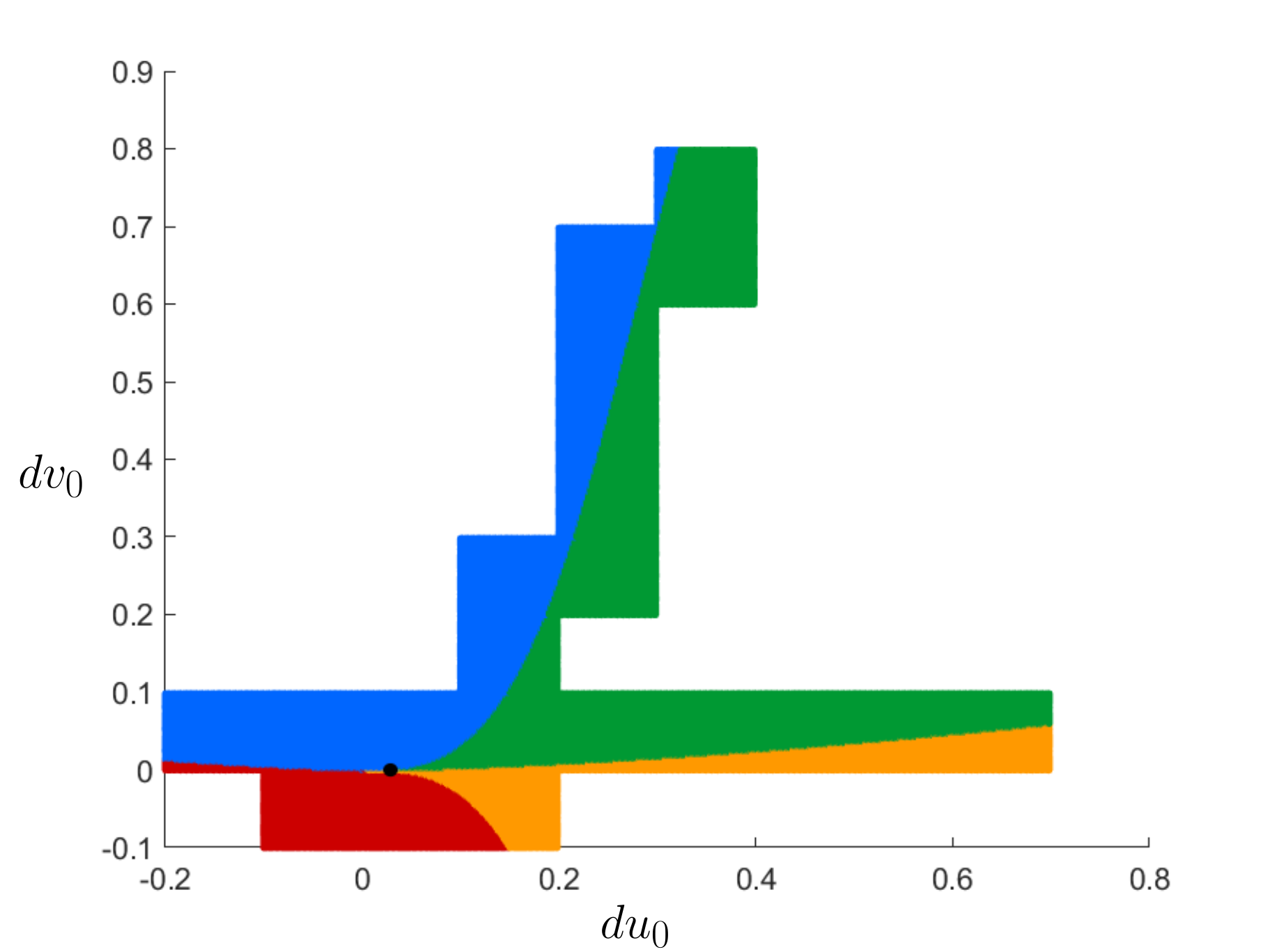}
			\caption{Diagram for $\lambda = 1$ (dense grid). The pair of parameters corresponding to a solution is marked by a black dot}
			\label{fig: shooting, lambda = 1, dense}
		\end{subfigure} \vspace{0.2 cm}
		\caption{The $du_0$-$dv_0$ diagram for $\lambda = 1$ and both coarse and dense grid}
		\label{fig: shooting, lambda = 1 LOWER}
	\end{figure}

	\begin{figure}[phtb!]
		\centering
		\begin{subfigure}{0.4\textwidth}
			\centering
			\includegraphics[width = \linewidth]{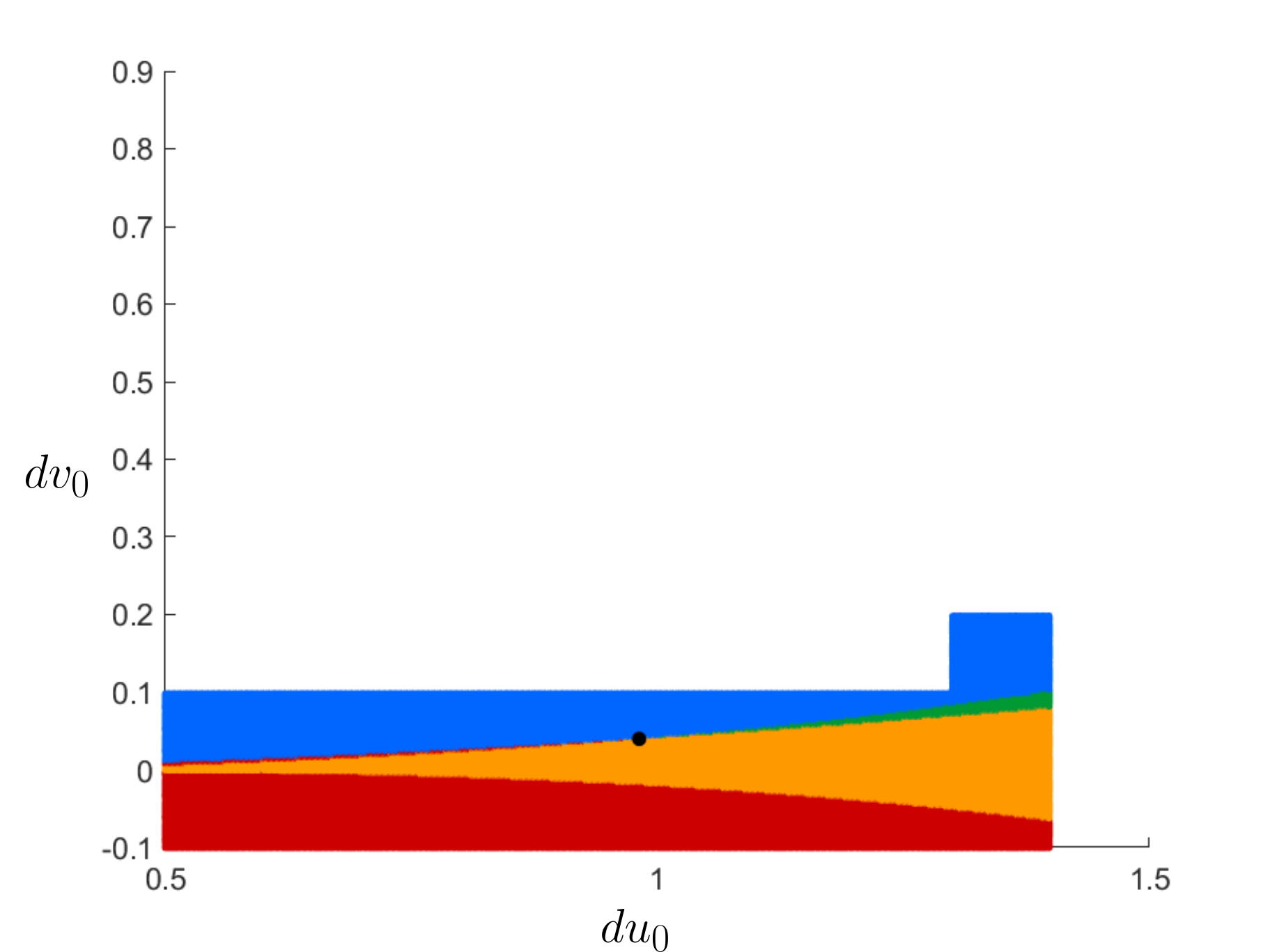}
			\caption{The $du_0$-$dv_0$ diagram for $\lambda = 10$,\\$du_0 \approx 1$, $dv_0 \approx 0.03$}
			\label{fig: shooting, lambda = 10, dense LOWER}
		\end{subfigure}
		\begin{subfigure}{0.4\textwidth}
			\centering
			\includegraphics[width = \linewidth, keepaspectratio]{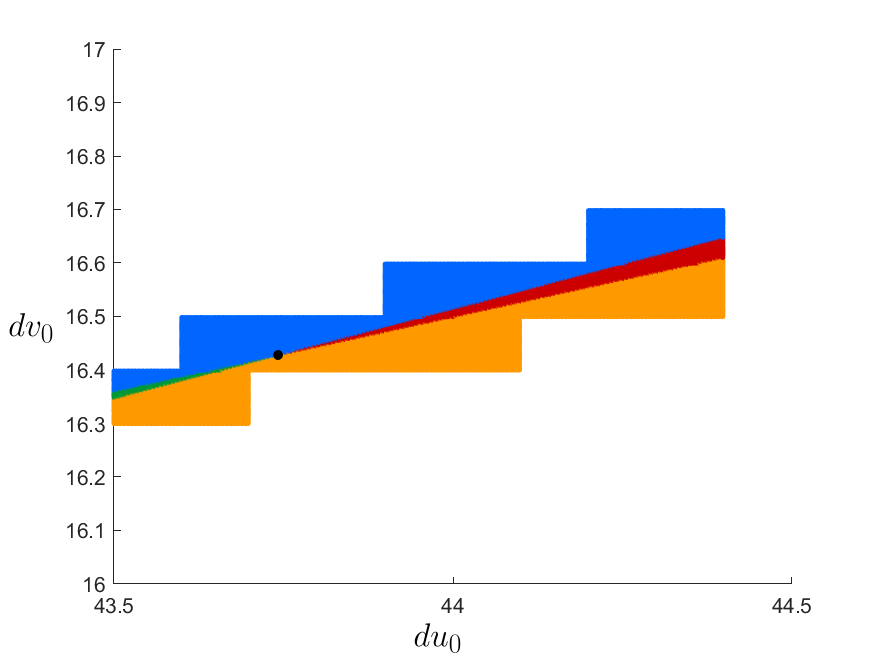}
			\caption{The $du_0$-$dv_0$ diagram for $\lambda = 10$,\\$du_0 \approx 44$, $dv_0 \approx 16.5$}
			\label{fig: shooting, lambda = 10, dense UPPER}
		\end{subfigure}
		\ \\
		\begin{subfigure}{0.4\textwidth}
			\centering
			\includegraphics[width = \linewidth]{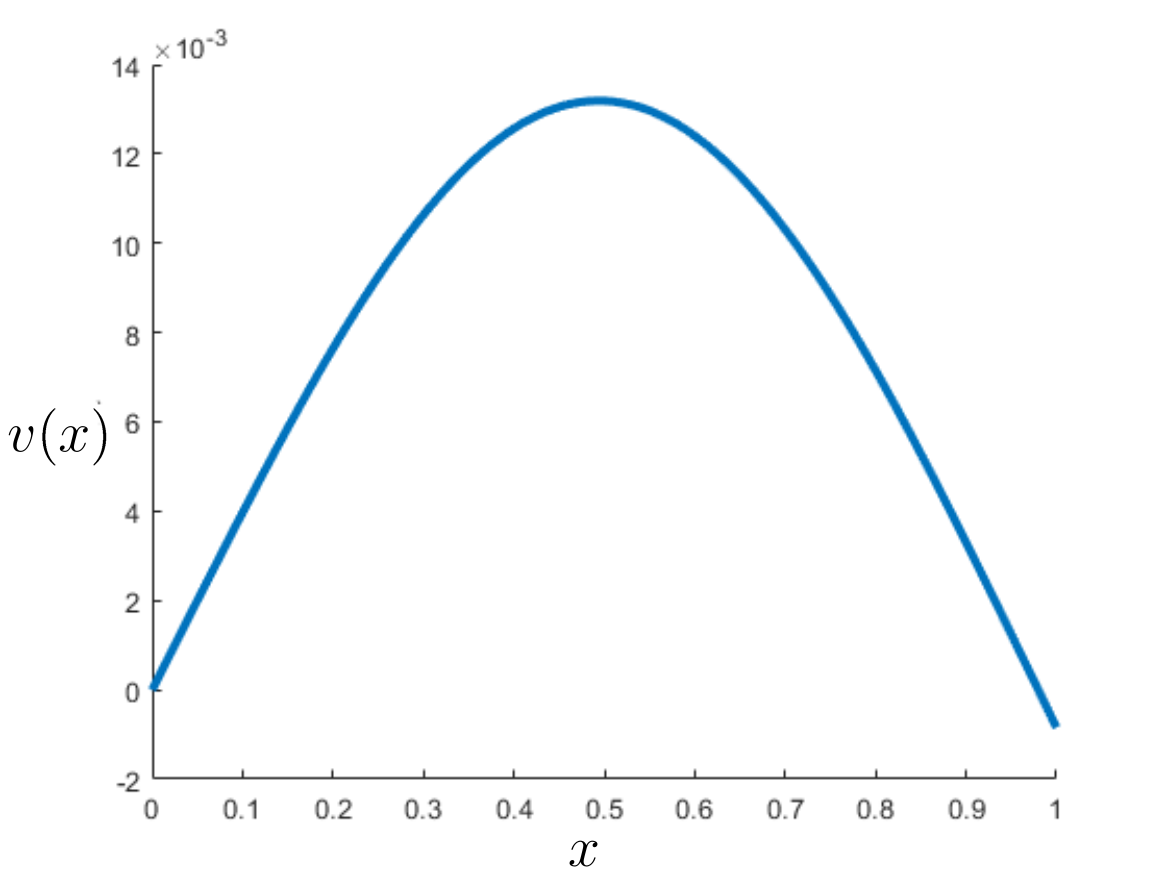}
			\caption{Lower branch solution for $\lambda = 10$\\corresponding to Figure \ref{fig: shooting, lambda = 10, dense LOWER}}
			\label{fig: lower solution, lambda = 10}
		\end{subfigure} 
		\begin{subfigure}{0.4\textwidth}
			\centering
			\includegraphics[width = \linewidth]{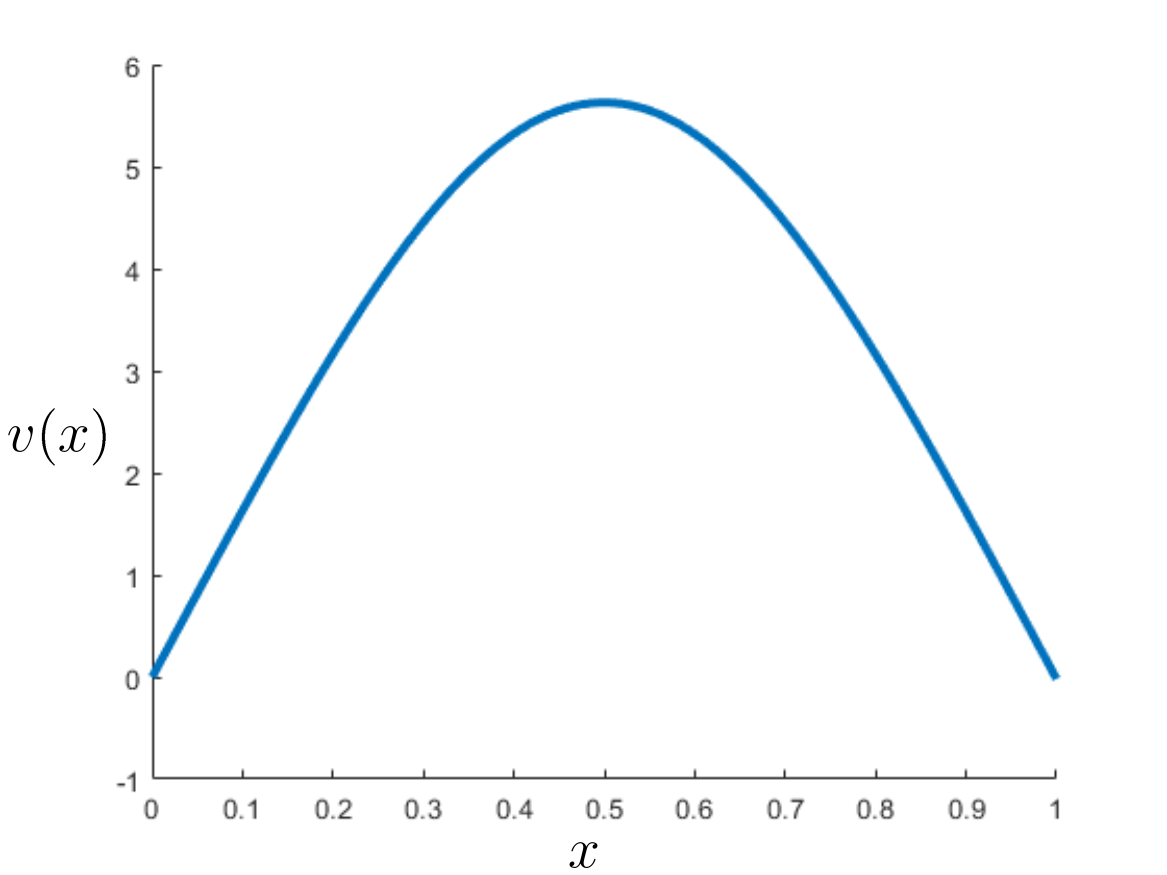}
			\caption{Upper branch solution for $\lambda = 10$\\corresponding to Figure \ref{fig: shooting, lambda = 10, dense UPPER}}
			\label{fig: upper solution, lambda = 10}
		\end{subfigure}
		\caption{Numerical experiments for $\lambda = 10$}
        \label{fig: shooting, lambda = 10}
	\end{figure}

	Figure \ref{fig: shooting, lambda = 1, dense} shows a narrow red protrusion coming from the red area at the bottom-left corner of the diagram. At the point where the protrusion touches the green area, all four colors connect at the pair $(du_0, dv_0)$ corresponding to a nontrivial numerical solution.

	As the value of $\lambda$ increases, the red protrusion is visible for higher and higher values of $du_0$ and $dv_0$, as shown
	in Figure \ref{fig: shooting, lambda = 10, dense LOWER}. 
	
	If we fix $\lambda = 10$, besides the solution illustrated in Figure \ref{fig: shooting, lambda = 10, dense LOWER} with $du_0 \approx 1$ and $dv_0 \approx 0.03$, 
	the experiments found another solution for $du_0 \approx 44$, $dv_0 \approx 16.5$, as illustrated in Figure \ref{fig: shooting, lambda = 10, dense UPPER}. Also, Figures~\ref{fig: lower solution, lambda = 10}, \ref{fig: upper solution, lambda = 10} show that the $L^{\infty}$-norm of the solution for lower initial slopes $du_0$, $dv_0$ (corresponding to the diagram in Figure 
	\ref{fig: shooting, lambda = 10, dense LOWER}) is lower than the $L^{\infty}$-norm of the other solution. The numerical experiments yield a similar scenario when $\lambda = 1$.\\
	
The numerical results above discussed strongly suggested the existence of two branches in the bifurcation diagram; the \emph{lower branch} (closer to the trivial solution -- in the view of the $L^{\infty}$-norm) 
	and the \emph{upper branch} (farther from the trivial solution).
	Based on the presented results, the branches are also getting closer to each other with increasing $\lambda$, presumably colliding when $\lambda=\lambda_{\rm bif}$. 
	
\begin{figure}[htb]
		\centering
		\begin{subfigure}{0.3\textwidth}
			\centering
			\includegraphics[width = \textwidth]{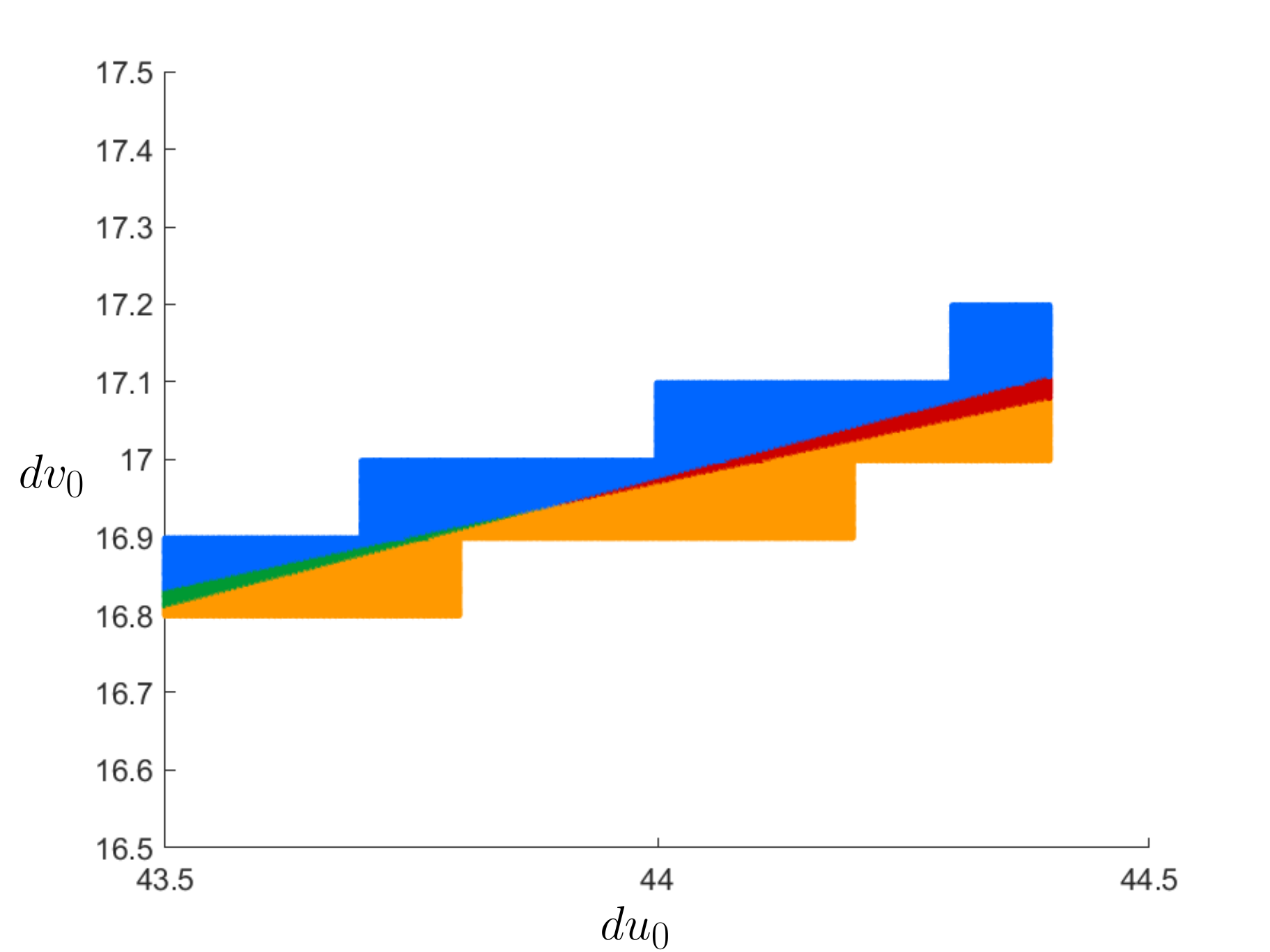}
			\caption{Solution $v^{5}$ in $du_0$-$dv_0$ diagram\\ \ }
			\label{fig: diagram upper lambda = 5}
		\end{subfigure}
		\begin{subfigure}{0.3\textwidth}
			\centering
			\includegraphics[width = \textwidth]{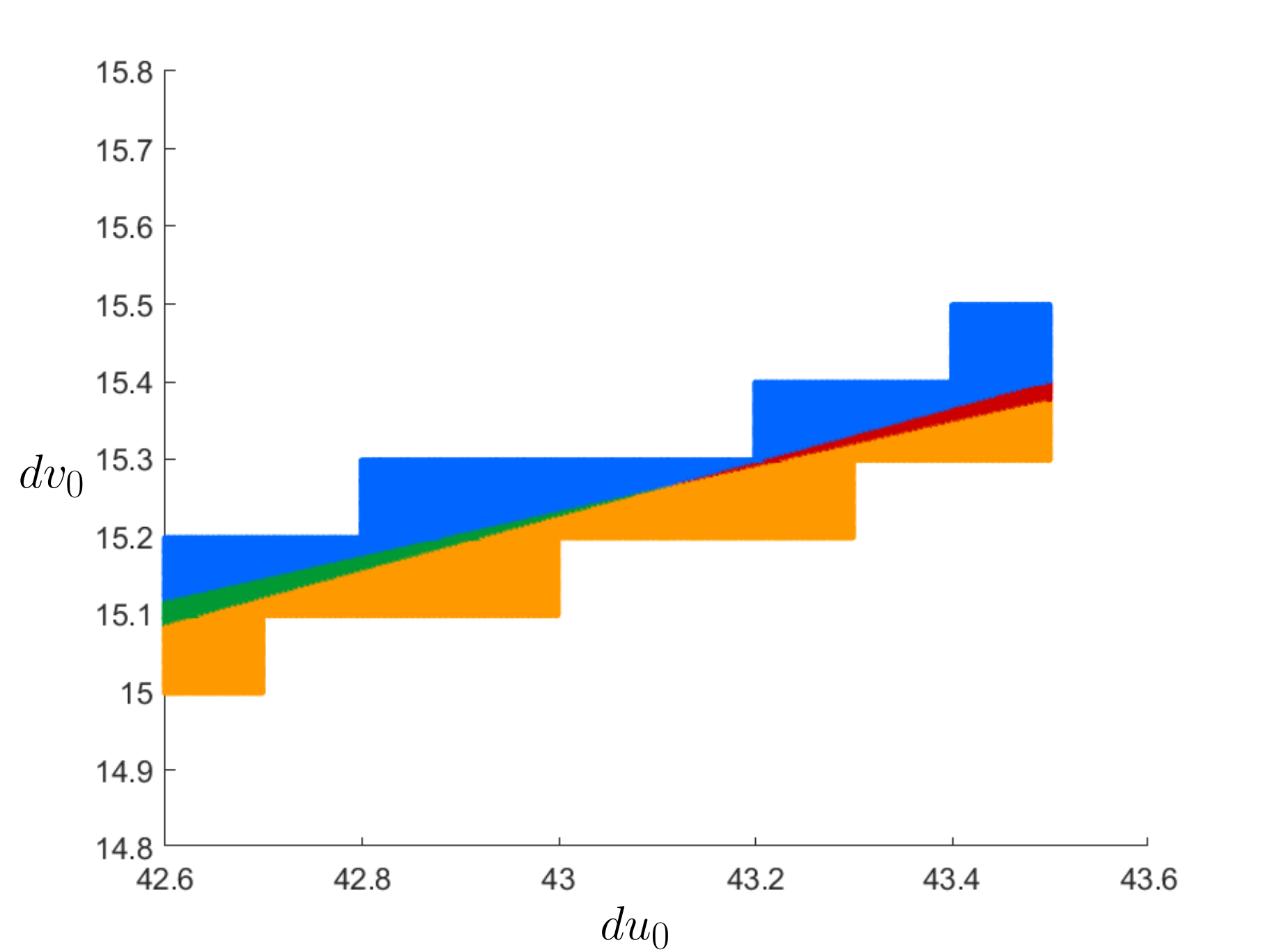}
			\caption{Solution $v^{20}$ in $du_0$-$dv_0$ diagram}
			\label{fig: diagram upper lambda = 20}
		\end{subfigure} 
		\begin{subfigure}{0.3\textwidth}
			\centering
			\includegraphics[width = \textwidth]{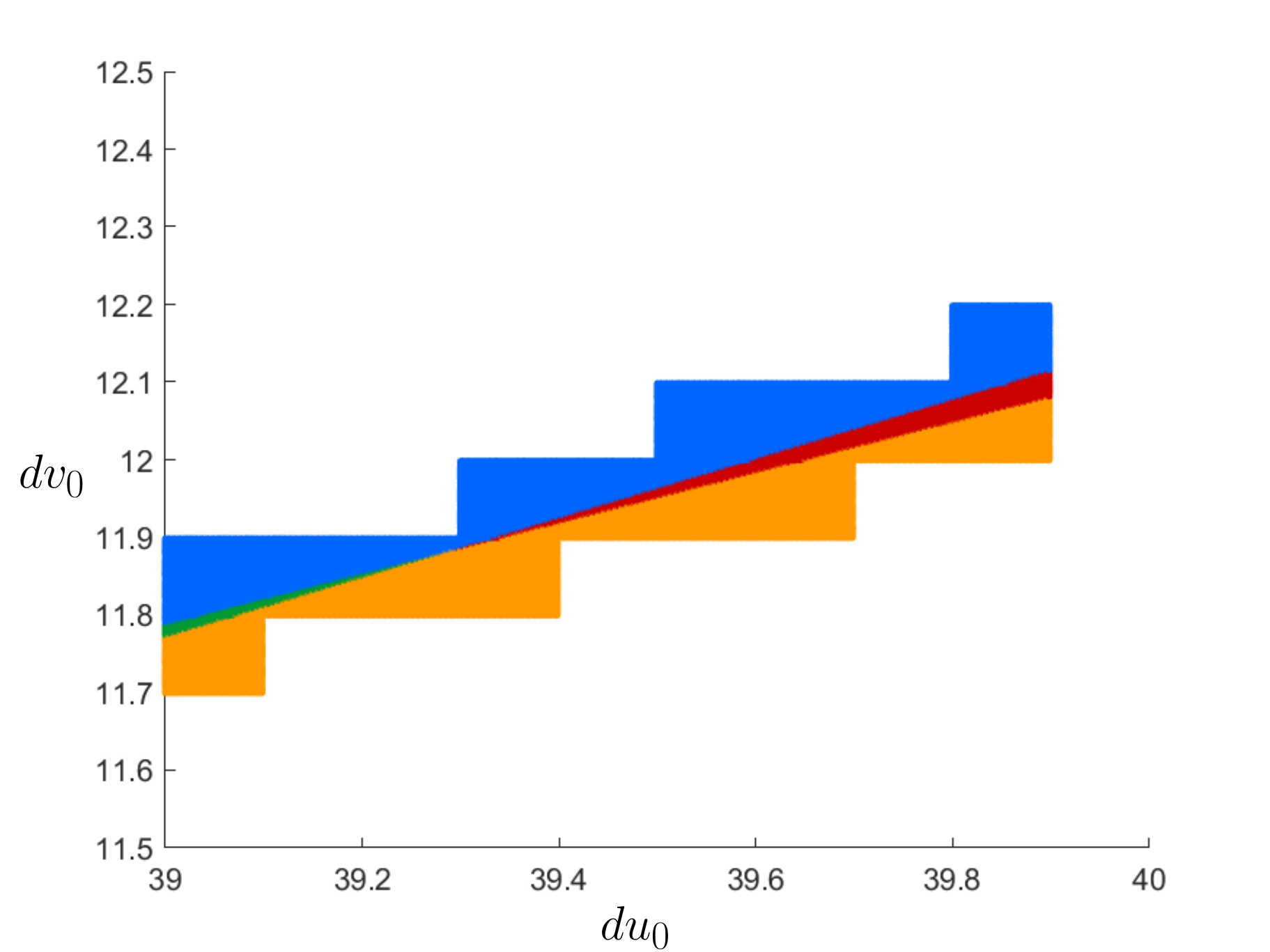}
			\caption{Solution $v^{40}$ in $du_0$-$dv_0$ diagram}
			\label{fig: diagram upper lambda = 40}
		\end{subfigure}
		\caption{Solutions from upper branch of the bifurcation diagram shown in $du_0$-$dv_0$ diagram}
		\label{fig: diagram upper branch}
	\end{figure}

\begin{figure}[htb]
		\centering
		\begin{subfigure}{0.3\textwidth}
			\centering  
			\includegraphics[width = \textwidth]{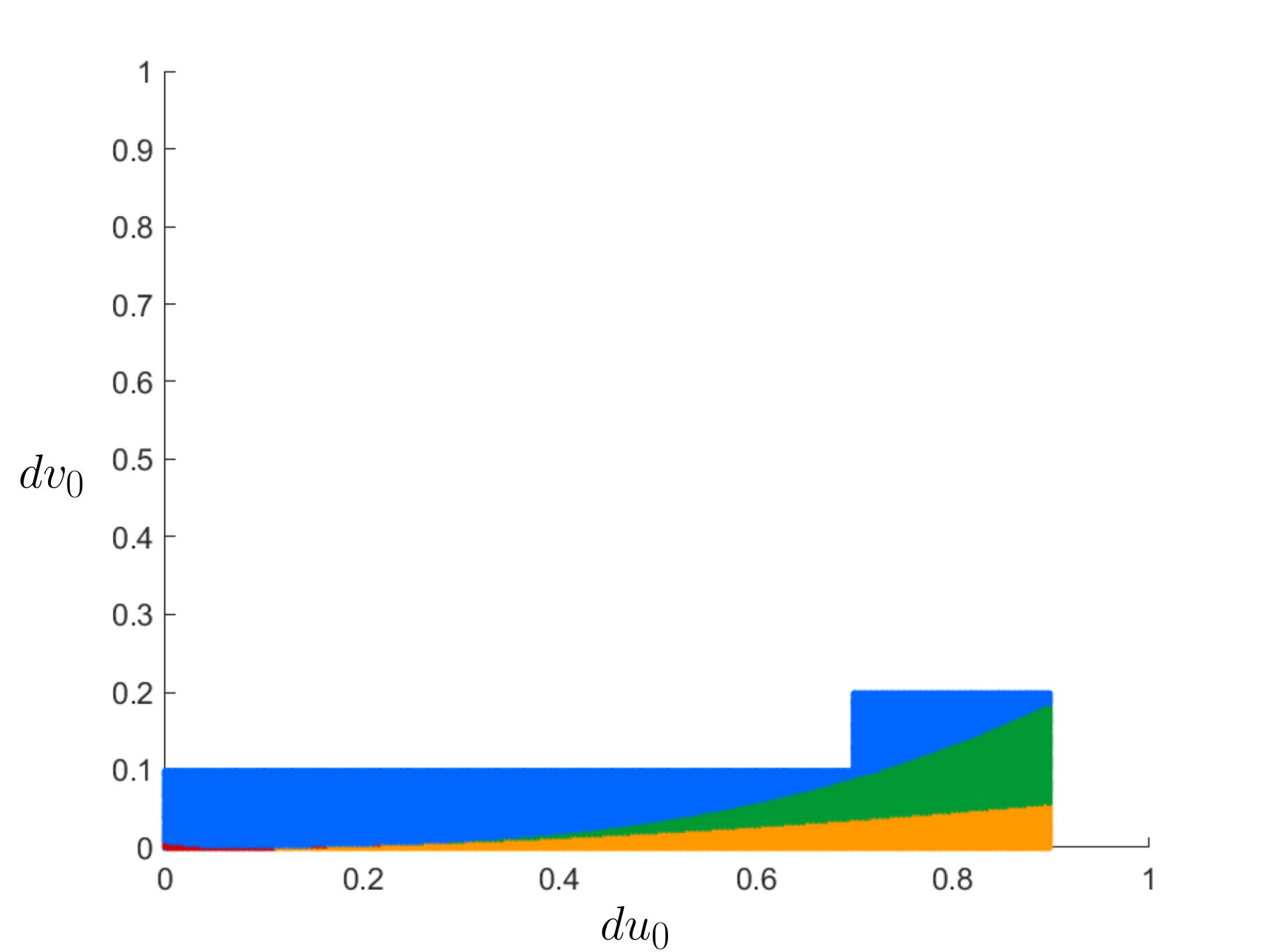}
			\caption{Solution $v_{5}$ in $du_0$-$dv_0$ diagram\\ \ }
			\label{fig: diagram lower lambda = 5}
		\end{subfigure}
		\begin{subfigure}{0.3\textwidth}
			\centering
			\includegraphics[width = \textwidth]{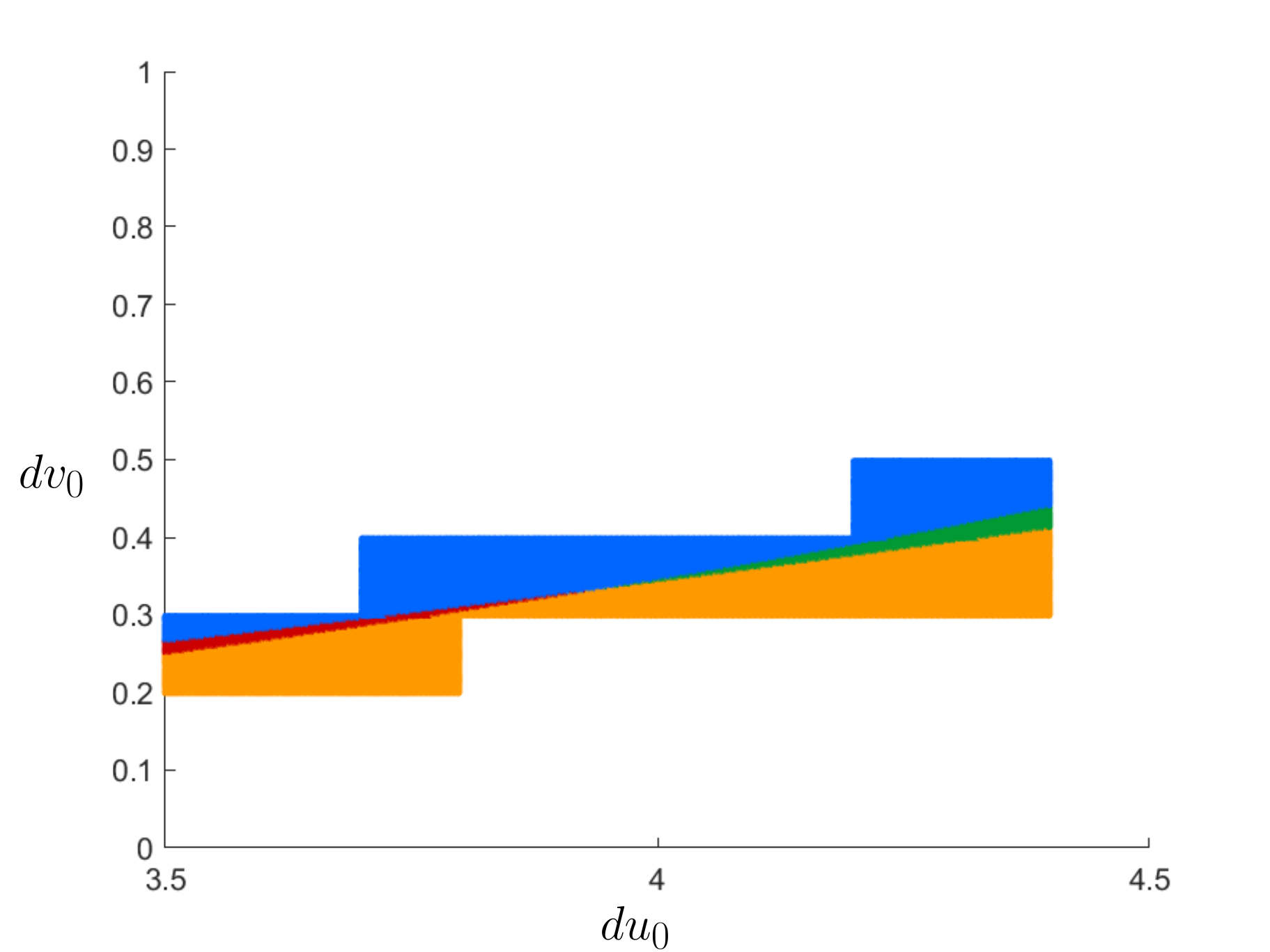}
			\caption{Solution $v_{20}$ in $du_0$-$dv_0$ diagram}
			\label{fig: diagram lower lambda = 20}
		\end{subfigure}
		\begin{subfigure}{0.3\textwidth}
			\centering
			\includegraphics[width = \textwidth]{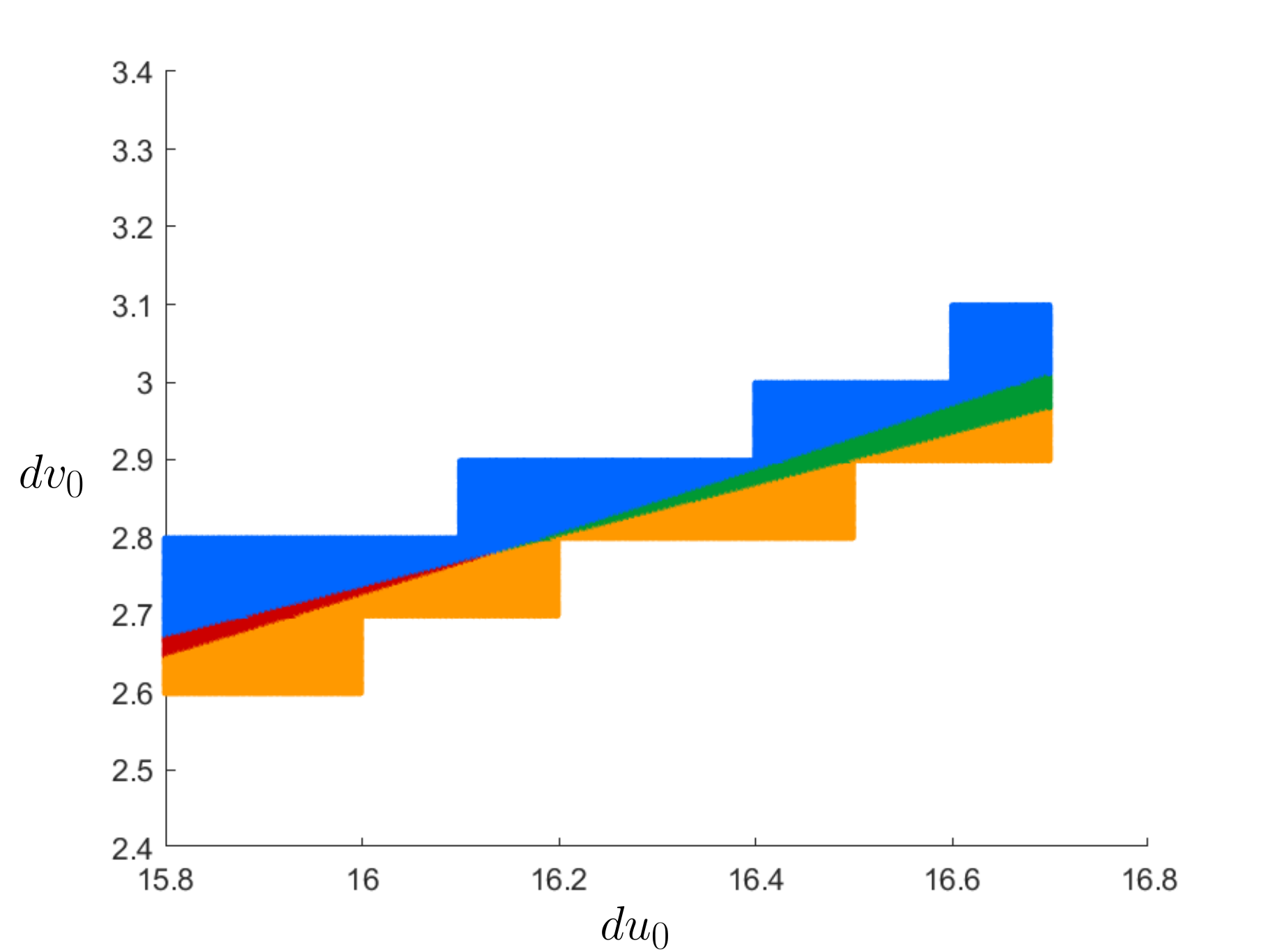}
			\caption{Solution $v_{40}$ in $du_0$-$dv_0$ diagram}
			\label{fig: diagram lower lambda = 40}
		\end{subfigure}
		\caption{Solutions from lower branch of the bifurcation diagram shown in $du_0$-$dv_0$ diagram}
		\label{fig: diagram lower branch}
	\end{figure}

	With this assumption, we let the script explore the range $1 \leq \lambda \leq 50$. The results confirmed the assumptions that the branches meet at some point
	$\lambda_{\rm bif} \approx 49$. The results also confirmed (numerically) that beyond $\lambda_{\rm bif}$, no solutions could be found. This behavior is summarized in the bifurcation diagram in Figure \ref{fig:bifurcation diagram}.
	Further illustrations can be found in Figures \ref{fig: diagram upper branch} and \ref{fig: diagram lower branch}.


	\FloatBarrier
\bigskip

{\bf Acknowledgements.} The authors were supported by the Grant CR 22-18261S of the Grant Agency of the Czech Republic.

\bigskip

{\bf Conflict of Interest.} The authors declare that they have no conflict of interest. 
   	\bibliographystyle{abbrvdin}
	\bibliography{reference}	
\end{document}